\documentclass[review,3p,square]{elsarticle}

\usepackage{amssymb}
\usepackage{amsmath}
\usepackage{amsthm}
\usepackage{stmaryrd}
\usepackage{bbm}
\usepackage{enumitem}
\usepackage{mathrsfs}

\usepackage{fancyhdr}
\usepackage{hyperref}
\hypersetup{%
  pdftitle={BSDEs},%
  pdfsubject={BSDEs},%
  pdfauthor={Yan Wang, Xinying Li, Chuang Gu, Shengjun Fan},%
  pdfkeywords={BSDEs},%
  pdfstartview=FitH,%
  CJKbookmarks=true,%
  bookmarksnumbered=true,%
  bookmarksopen=true,%
  colorlinks=true, linkcolor=blue, urlcolor=blue, citecolor=blue, %
}

\usepackage{cleveref}
\crefname{thm}{Theorem}{Theorems}
\crefname{pro}{Proposition}{Propositions}
\crefname{lem}{Lemma}{Lemmas}
\crefname{rmk}{Remark}{Remarks}
\crefname{cor}{Corollary}{Corollaries}
\crefname{dfn}{Definition}{Definitions}
\crefname{ex}{Example}{Examples}
\crefname{section}{Section}{Sections}
\crefname{subsection}{Subsection}{Subsections}

\newcommand{\To}{\rightarrow}
\newcommand{\as}{{\rm d}\mathbb{P}\times{\rm d} t-a.e.}

\newcommand{\dif}{{\rm d}}
\newcommand{\ps}{\mathbb{P}-a.s.}

\newcommand{\essinf}{\mathop{\operatorname{ess\,inf}}}
\newcommand{\F}{\mathcal{F}}
\newcommand{\E}{\mathbb{E}}

\newcommand{\M}{\mathcal{M}}

\newcommand{\s}{\mathcal{S}}

\newcommand{\T}{[0,T]}

\newcommand{\A}{\mathcal{A}}
\newcommand{\R}{{\mathbb R}}
\newcommand{\Q}{{\mathbb Q}}

\newcommand{\sgn}{{\rm{sgn}}}

\newtheorem{theorem}{Theorem}[section]

\newtheorem{lemma}[theorem]{Lemma}

\newtheorem{remark}[theorem]{Remark}
\newtheorem{proposition}[theorem]{Proposition}
\newtheorem{example}[theorem]{Example}

\theoremstyle{definition}

\journal{arXiv}

\begin{document}
\begin{frontmatter}

\title{{On the existence and uniqueness of unbounded solutions \\ to quadratic BSDEs with convex generators}\tnoteref{found}}
\tnotetext[found]{This work is supported by the National Natural Science Foundation of China (Grant No. 12171471).
\vspace{0.2cm}}


\author[a,b]{Yan Wang} \ead{wangyan\_shuxue@163.com}
\author[b]{Xinying Li} \ead{lixinyingcumt@163.com}
\author[b]{Chuang Gu} \ead{guchuang1026@163.com}
\author[b]{Shengjun Fan\corref{cor}} \ead{shengjunfan@cumt.edu.cn} \vspace{-0.5cm}

\affiliation[a]{organization={Jiangsu University of Architectural Technology} ,
            city={Xuzhou 221116},
            country={China}}

\affiliation[b]{organization={School of Mathematics, China University of Mining and Technology},
            city={Xuzhou 221116},
            country={China}}

\cortext[cor]{Corresponding author\vspace{0.2cm}}

\begin{abstract}
With the terminal value $\xi^-$ admitting a certain exponential moment and $\xi^+$ admitting all exponential moments or being bounded, we establish several existence and uniqueness results for unbounded solutions of backward stochastic differential equations (BSDEs) whose generator $g$ satisfies a monotonicity condition with a general growth in the first unknown variable $y$ and a convexity condition with a quadratic growth in the second unknown variable $z$. In particular, the generator $g$ may not be locally-Lipschitz continuous in $y$. This strengthens some results reported in Delbaen, Hu and Richou [Ann. Inst. Henri Poincar$\acute{e}$, 2011(47), 559-574] by relaxing the continuity and growth of $g$ in $y$. We also give an explicit expression of the first process in the unique unbounded solution of a BSDE when the generator $g$ is jointly convex in $(y,z)$ and has a linear growth in $y$ and a quadratic growth in $z$. Finally, we prove the corresponding comparison theorems for unbounded solutions of the preceding BSDEs.
 \vspace{0.2cm}
\end{abstract}

\begin{keyword}
Backward stochastic differential equation \sep Existence and uniqueness \sep  Unbounded solution \sep \hspace*{1.75cm} Quadratic growth \sep Convex generator.
\vspace{0.2cm}

\MSC[2010] 60H10\vspace{0.2cm}
\end{keyword}

\end{frontmatter}
\vspace{-0.4cm}

\section{Introduction} \label{section1}
Throughout this paper, let us fix a positive integer $d$ and a positive real $T$. Let $(\Omega,\F,\mathbb{P})$ be a complete probability space carrying a standard $d$-dimensional Brownian motion $(B_t)_{t\geq0}, (\F_t)_{t\geq0}$ the completed natural $\sigma$-algebra generated by $(B_t)_{t\geq0}$ and $\F_T := \F$. We consider the following one-dimensional quadratic backward stochastic differential equation (BSDE in short):
\begin{align}\label{BSDE1.1}
     Y_t=\xi-\int_t^Tg(s,Y_s,Z_s)\dif s+\int_t^TZ_s\cdot\dif B_s,~~~~t\in[0,T],
\end{align}
where $T$ is the terminal time, and the terminal value $\xi$ is an $\F_T$-measurable unbounded random variable. The generator
$$
g(\omega,t,y,z):\Omega\times[0,T]\times\R\times\R^d\mapsto\R
$$
is $(\F_t)$-progressively measurable for each $(y,z)\in\R\times\R^d$, which is continuous in $(y,z)$ and has a quadratic growth in $z$. The BSDE with parameters $(\xi,T,g)$ is denoted by BSDE $(\xi,T,g)$.

In 1990, the seminal paper \cite{Pardoux-Peng 1990} first introduced the nonlinear BSDEs and established an existence and uniqueness result for the adapted solution of BSDEs under the uniformly Lipschitz continuous condition of the generator, which laid the foundation for the theoretical and applied research of BSDEs. After that, a lot of efforts have been made in order to study the well-posedness of these equations and their practical applications. Since it appears naturally in the utility maximization problem, quadratic BSDE is a class of BSDEs that have attracted particular attention; see for example \cite{Hu 2005, Rouge 2000, Fan 2016, Luo-Fan 2018} along with the updated papers \cite{Tian 2023, Hu-Tang-Wang 2022, Fan-Hu-Tang 2023} for more details. In particular, Kobylanski \cite{Kobylanski 2000} first studied quadratic BSDEs with bounded terminal values, and established a rather general existence and uniqueness result on the bounded solution. This result was further extended in \cite{Fan 2016} to quadratic BSDEs with finite and infinite horizons, where the generator $g$ can satisfy a monotonicity condition with a general growth in $y$ and a locally-Lipschitz continuous condition with a quadratic growth in $z$. Roughly speaking, the monotonicity condition of the generator $g$ in $y$ can be read as follows: there exists a constant $\beta\geq 0$ such that $\as$, for each $(y_1,y_2,y,z)\in\R\times\R\times\R\times\R^d$,
\begin{equation}\label{eq:1.2}
-\sgn(y_1-y_2)\big(g(t,y_1,z)-g(t,y_2,z)\big)\leq\beta|y_1-y_2|.
\end{equation}

By applying the localization method, the authors in \cite{Briand-Hu 2006} obtained the first existence result for the unbounded solution of a quadratic BSDE with unbounded terminal value admitting a certain exponential moment. Based on the $\theta$-difference technique for convex functions, the uniqueness was established in \cite{Briand-Hu 2008} when the terminal value possesses all exponential moments and the generator $g$ is uniformly Lipschitz continuous in $y$ and convex in $z$. By virtue of the Legendre-Fenchel transform of convex functions, this result was further strengthened in \cite{Delbaen 2011}, which proved that the uniqueness holds among some proper space when the terminal value admits a certain exponential moment. Furthermore, the critical case was tackled in \cite{Delbaen 2015}, but an additional assumption that the generator $g$ is independent of $y$ and strongly convex in $z$ was required. We would like to especially mention that in \cite{Briand-Hu 2008} and \cite{Delbaen 2011}, the generator $g$ is forced to satisfy the uniformly Lipschitz continuous condition and then the linear growth condition in the first unknown variable $y$ in order to guarantee uniqueness of the unbounded solution.

Recently, the authors in \cite{Fan-Hu-Tang 2020} proved a uniqueness result for the unbounded solution of a quadratic BSDE under the assumptions that the terminal value $\xi$ has all exponential moments and either $f(t,y,z):=-g(t,y,-z)$ or $f(t,y,z):=g(t,-y,z)$ satisfies the following extended convexity condition in the state variables $(y,z)$: there exist two constants $C\geq 0$ and $\delta\in (0,1)$ such that for each $\theta\in (0,1)$ and each $(t,y_1,y_2,z_1,z_2)\in \T\times\R\times\R\times\R^d\times\R^d$, we have
\begin{align}\label{eq:1.3}
&\mathbbm{1}_{\{y_1-\theta y_2\geq 0\}}\frac{f(t,y_1,z_1)-\theta f(t,y_2,z_2)}{1-\theta}\nonumber\\
&\ \ \ \leq C\left(1+|y_1|+|y_2|+|z_1|^{1+\delta}+|z_2|^{1+\delta}+\left|\frac{y_1-\theta y_2}{1-\theta}\right|+\left|\frac{z_1-\theta z_2}{1-\theta}\right|^2\right),
\end{align}
which admits a general growth of $g$ in $y$ and weakens, at the first time, the uniformly Lipschitz continuous condition of $g$ in $y$ used in \cite{Briand-Hu 2008,Delbaen 2011}, to the best of our knowledge. Note from \cite[Remark 8]{Fan-Hu-Tang 2020} that if \eqref{eq:1.3} is fulfilled, then the generator $g$ has to satisfy the monotonicity condition \eqref{eq:1.2}. On the other hand, it is not very difficult to verify that the following generator
\begin{equation}\label{eq:1.4}
g(t,y,z):=\mathbbm{1}_{y<1} (y-1)^3+\mathbbm{1}_{y\geq 1}\sqrt{y-1}+\frac{1}{2}|z|^2,\ \ \ (t,y,z)\in\T\times\R\times\R^d
\end{equation}
satisfies the monotonicity condition \eqref{eq:1.2} with a general growth, but neither $f(t,y,z):=-g(t,y,-z)$ nor $f(t,y,z):=g(t,-y,z)$ satisfies \eqref{eq:1.3}. Indeed, since the function $g(t,\cdot,z)$ is nondecreasing on $\R$ for each $(t,z)\in \T\times\R^d$, then \eqref{eq:1.2} with $\beta=0$ follows immediately. In what follows, we verify the assertion that neither $f(t,y,z):=-g(t,y,-z)$ nor $f(t,y,z):=g(t,-y,z)$ satisfies \eqref{eq:1.3}. In fact, by letting $y_1=y_2=1$, $z_1=\frac{n-1}{n}z$, $z_2=z$ and $\theta=\frac{n-1}{n}$ for an integer $n\geq 1$ and a vector $z\in\R^d$ with $|z|=n$, \eqref{eq:1.3} becomes
\begin{align}\label{eq:1.5}
&\mathbbm{1}_{\{1-\theta\geq 0\}}\frac{-g(t,1,-z_1)+\theta g(t,1,-z_2)}{1-\theta}=\frac{n}{2}
\bigg[-\big|\frac{n-1}{n}z\big|^2+\frac{n-1}{n}|z|^2\bigg]=\frac{n(n-1)}{2}\nonumber\\
&\ \ \ \leq C\big(4+|z_1|^{1+\delta}+|z_2|^{1+\delta}\big)=C\big(4+(n-1)^{1+\delta}
+n^{1+\delta}\big),
\end{align}
which is impossible when $n\To \infty$ due to the fact $1+\delta<2$. Consequently, $f(t,y,z):=-g(t,y,-z)$ fails to satisfy \eqref{eq:1.3} due to the presence of the term $\frac{1}{2}|z|^2$. Furthermore, by letting $y_1=n-1$, $y_2=n$ and $\theta=\frac{n-1}{n}$ for an integer $n\geq 1$, \eqref{eq:1.3} becomes\vspace{0.2cm}
\begin{align*}
&\mathbbm{1}_{\{y_1-\theta y_2\geq 0\}}\frac{g(t,-y_1,0)-\theta g(t,-y_2,0)}{1-\theta}=n\bigg[(-n)^3-\frac{n-1}{n}(-n-1)^3\bigg]=2n^3-2n-1\\
&\ \ \ \leq C(1+n-1+n\big)=2nC,
\end{align*}
which is impossible when $n\To \infty$. Consequently, $f(t,y,z):=g(t,-y,z)$ fails to satisfy \eqref{eq:1.3} due to the presence of the term $\mathbbm{1}_{y<1}(y-1)^3$. Thus, the desired assertion is verified. In addition, we would like to emphasize that the generator $g$ defined in \eqref{eq:1.4} does not satisfy the locally-Lipschitz continuity condition in $y$ due to the presence of the term $\mathbbm{1}_{y\geq 1}\sqrt{y-1}$.


It can be concluded from the last paragraph that the monotonicity condition \eqref{eq:1.2} along with a general growth is strictly weaker than the extended convexity condition \eqref{eq:1.3} proposed in \cite{Fan-Hu-Tang 2020}. Then, in light of \cite{Fan 2016} and \cite{Fan-Hu-Tang 2020}, a question is naturally asked: under the elementary assumptions that the generator $g$ satisfies the monotonicity condition \eqref{eq:1.2} along with a general growth in $y$ and the convexity condition with a quadratic growth in $z$, without imposing the locally-Lipschitz continuity condition of $g$ in $y$ nor the extended convexity condition \eqref{eq:1.3}, does the uniqueness for the unbounded solution of a quadratic BSDE hold still? The present paper gives some affirmative answers. More specifically, the main purpose of this paper is to study the existence and uniqueness for unbounded solutions of quadratic BSDEs when the terminal value $\xi^-$ admits a certain exponential moment, and the generator $g$ may not be locally-Lipschitz continuous in $y$, and is convex in $z$. Our result strengthens those established in \cite{Delbaen 2011}, by relaxing the continuity and growth condition of the generator $g$ in the state variable $y$.

The paper is organized as follows. In Proposition \ref{cun1} of Section \ref{section 2}, we prove a general existence result for the unbounded solution of BSDE (\ref{BSDE1.1}) in the spirit of \cite{Fan-Hu-Tang 2020}, which can be compared with Corollary 2.2 in \cite{Delbaen 2011}. Here we work with $\xi^-$ admitting a certain exponential moment and $\xi^+$ admitting all exponential moments or being bounded. Section \ref{section 3} is devoted to the uniqueness for the unbounded solution of BSDE (\ref{BSDE1.1}) in the case of the generator $g$ satisfying the monotonicity condition \eqref{eq:1.2} with general growth in $y$ and the convexity condition with a quadratic growth in $z$. The generator $g$ defined in \eqref{eq:1.4} is a typical example satisfying all these conditions, without imposing the locally-Lipschitz continuity condition of $g$ in $y$ nor the extended convexity condition \eqref{eq:1.3}. The main result is stated in Theorem \ref{unth3.1}. We remark that the authors in \cite{Delbaen 2011} proved the uniqueness result when both $\xi^+$ and $\xi^-$ have a certain exponential moment and the generator $g$ is uniformly Lipschitz continuous in $y$ and convex with a quadratic growth in $z$. In contrast, we relax the continuity and growth condition of $g$ in $y$, and lift the requirement for the integrability of $\xi^+$, but actually this requirement is not necessary; see (ii) of Remark \ref{rmk-hfsj} for more details. In Example \ref{lz}, we provide several examples to which Theorem \ref{unth3.1} but no existing results can apply. In Theorem \ref{unthc} of Section \ref{section 4}, we prove the uniqueness for the unbounded solution of a quadratic BSDE in the case of $g$ being jointly convex in $(y,z)$. Actually, the assumptions used in Theorem \ref{unthc} are strictly stronger than those in \cite{Delbaen 2011}. However, we get a more accurate conclusion in this setting, that is, $Y_\cdot$ admits a more explicit expression; see (ii) of Remark \ref{rmk-lip} for more details. In order to prove Theorems \ref{unth3.1} and \ref{unthc}, we systematically apply some techniques and ideas used in \cite{Delbaen 2011} and \cite{Delbaen 2015}, including the Legendre-Fenchel transform of convex functions, Girsanov's theorem, the de La Vall$\acute{\rm{e}}$e Poussin lemma and Fenchel's inequality. Meanwhile, we adopt an exponential shifting transform to deal with the monotonicity condition of the generator $g$ in $y$ as well as the general growth condition. It should be especially emphasized that although most of these ideas and techniques have been used in BSDE literature for a long time, combining them together in a novel way to prove the desired assertions under our weaker assumptions on the coefficients of BSDEs is also a contribution of this paper. Finally, in Theorem \ref{cth} of Section \ref{section 5}, it seems to be the first time, to the best of our knowledge, that the corresponding comparison theorems for unbounded solutions of quadratic BSDEs are established under the assumptions used in Sections \ref{section 3} and \ref{section 4}, with the help of Proposition 5 proposed in \cite{Briand-Hu 2006}, the uniqueness stated in Theorems \ref{unth3.1} and \ref{unthc}, and the infinite evolution technique combined with the localization argument. The conclusion of this paper is presented in Section \ref{section 6}.

Let us close this introduction with some notations that will be used later. Let $x\cdot y$ denote the Euclidean inner product for $x, y\in\R^d$. Denote by $\mathbbm{1}_A$ the indicator of set $A$, and $\sgn(x):=\mathbbm{1}_{x>0}-\mathbbm{1}_{x\leq0}$. Let $a\wedge b$ be the minimum of $a$ and $b$, $a^-:=-(a\wedge0)$ and $a^+:=(-a)^-$. Moreover, let $\R_+:=[0,+\infty)$, and $\R_-:=(-\infty,0]$. By an adapted solution to BSDE (\ref{BSDE1.1}), we mean a pair of $(\F_t)$-progressively measurable processes $(Y_t,Z_t)_{t\in[0,T]}$ valued in $\R\times\R^d$ such that $\ps$, the function $t\rightarrow Y_t$ is continuous, $t\rightarrow Z_t$ is square-integrable, $t\rightarrow g(t,Y_t,Z_t)$ is integrable, and $(Y_\cdot,Z_\cdot)$ verifies (\ref{BSDE1.1}). We denote by $\partial f$ the subdifferential of a function $f:\R^{d}\rightarrow\R$, and the subdifferential of $f$ at $z_0$ is the non-empty convex compact set of elements $u\in\R^{1\times d}$ such that
$$
f(z)-f(z_0)\geq u(z-z_0),~~ z\in\R^d.
$$
For any $(\F_t)$-progressively measurable process $(q_t)_{t\in[0,T]}$ such that $\ps, \int_0^T |q_s|^2 \dif s<+\infty$, we denote by $\mathcal{E}(q)$ the Dol{\rm{$\acute{e}$}}ans-Dade exponential
$$
\Big(\exp\Big(\int_0^tq_s\cdot\dif B_s-\frac{1}{2}\int_0^t |q_s|^2 \dif s\Big)\Big)_{t\in[0,T]}.
$$
Let us recall the following Fenchel's inequality
$$
xy\leq{\rm{exp}}(x)+y(\ln y-1),~~~ (x,y)\in\R\times\R^+.
$$
Then we have
\begin{align*}
xy=px\frac{y}{p}\leq{\rm{exp}}(px)+\frac{y}{p}(\ln y-\ln p-1).
\end{align*}
Furthermore, we define the following two process spaces.

$\bullet$ $\s^p([0,T];\R)$ is the set of all $(\F_t)$-progressively measurable and continuous real-valued processes $(Y_t)_{t\in[0,T]}$ satisfying
$$\|Y\|_{{\s}^p}:=\left(\E\left[\sup_{t\in[0,T]} |Y_t|^p\right]\right)^{1/p}<+\infty.\vspace{0.1cm}$$
We set $\s=\bigcup_{p>1}\s^p$.

$\bullet$ $\M^p([0,T];\R^d)$ is the set of all $(\F_t)$-progressively measurable $\R^d$-valued processes $(Z_t)_{t\in[0,T]}$ satisfying
$$
\|Z\|_{\M^p}:=\left\{\E\left[\left(\int_0^T |Z_t|^2{\rm d}t\right)^{p/2}\right] \right\}^{1/p}<+\infty.\vspace{0.3cm}
$$
Recall that an $(\F_t)$-progressively measurable real-valued process $(Y_t)_{t\in[0,T]}$ belongs to class $(D)$ if the $\{Y_\tau : \tau\in\Sigma_T\}$ is a family of uniformly integrable random variables, where $\Sigma_T$ represents the set of all $(\F_t)$-stopping times $\tau$ valued in $[0,T]$.

In the whole paper, we are always given two $(\F_t)$-progressively measurable $\R_+$-valued processes $(\underline{\alpha}_t)_{t\in[0,T]}$ and $(\overline{\alpha}_t)_{t\in[0,T]}$ satisfying
$$
\ps,\ \int_0^T(\underline{\alpha}_t+ \overline{\alpha}_t) \dif t<+\infty
$$
and two constants $\beta\geq0$ and $\gamma>0$ along with a continuous nondecreasing function $\varphi(\cdot):\R_+\rightarrow\R_+$ with $\varphi(0)=0$.

\section{An existence result}\label{section 2}
\setcounter{equation}{0}

This section will present and prove a general existence result for unbounded solutions of BSDE \eqref{BSDE1.1}. We first introduce the following assumption on the generator $g$.\vspace{0.2cm}

\noindent\textbf{(EX)} $\as$, for each $(y,z)\in\R\times\R^d$, we have
\begin{enumerate}
\renewcommand{\theenumi}{(1)}
\renewcommand{\labelenumi}{\theenumi}
\item $-\mathbbm{1}_{y>0}g(t,y,z)\leq \underline{\alpha}_t+\beta|y|+\frac{\gamma}{2}|z|^2$;
\renewcommand{\theenumi}{(2)}
\renewcommand{\labelenumi}{\theenumi}
\item $\mathbbm{1}_{y\leq0}g(t,y,z)\leq \overline{\alpha}_t+\beta|y|+\frac{\gamma}{2}|z|^2$;
\renewcommand{\theenumi}{(3)}
\renewcommand{\labelenumi}{\theenumi}
\item $|g(t,y,z)|\leq \overline{\alpha}_t+\varphi(|y|)+\frac{\gamma}{2}|z|^2$.
\end{enumerate}
\vspace{0.2cm}

The following existence proposition is the main result of this section, which can be compared with Theorem 2.1 and Corollary 2.2 in \cite{Delbaen 2011}.
\begin{proposition}\label{cun1}
Suppose that $\xi$ is a terminal value, the generator $g$ is continuous in $(y,z)$ and satisfies assumption {\rm{(EX)}}, and $\E[\exp(pe^{\beta T}(\xi^-+\int_0^T\overline{\alpha}_t\dif t))]<+\infty$ for some $p>\gamma$.
\begin{enumerate}[leftmargin=0.8cm]
\renewcommand{\theenumi}{{\rm{(i)}}}
\renewcommand{\labelenumi}{\theenumi}
\item If $\E[\exp(\overline{p}e^{\beta T}(\xi^++\int_0^T\underline{\alpha}_t\dif t))]<+\infty$ for some $\overline{p}>\gamma$, then BSDE {\rm{(\ref{BSDE1.1})}} admits a solution $(Y_\cdot,Z_\cdot)$ such that
\begin{align}\label{p1}
\E\Big[\exp\Big(p\sup_{t\in[0,T]}A_\beta(t)\Big)\Big]<+\infty,~ \text{where}~A_\beta(t):=e^{\beta t}\Big(Y_t^-+\int_0^t\overline{\alpha}_s \dif s\Big), ~t\in[0,T],
\end{align}
\vspace{-0.6cm}
\begin{align}\label{p2}
\E\Big[\exp\Big(\overline{p}\sup_{t\in[0,T]}\overline{A}_\beta(t)\Big)\Big]<+\infty, ~\text{where}~\overline{A}_\beta(t):=e^{\beta t}\Big(Y_t^++\int_0^t\underline{\alpha}_s \dif s\Big), ~t\in[0,T],
\end{align}
and $Z_\cdot\in\M^2$.
\renewcommand{\theenumi}{{\rm{(ii)}}}
\renewcommand{\labelenumi}{\theenumi}
\item If it holds that $\xi^+ +\int_0^T \underline{\alpha}_t \dif t\leq M$ for a constant $M>0$, then BSDE {\rm{($\ref{BSDE1.1}$)}} admits a solution $(Y_\cdot,Z_\cdot)$ satisfying {\rm{(\ref{p1})}} and $Y_\cdot \leq K$ for a constant $K>0$. Moreover, $Z_\cdot\in\M^2$.
\end{enumerate}
\end{proposition}
\noindent\textbf{Proof.}
We define, for any non-negative integrable function $l(\cdot):[0,T]\rightarrow\R_+$ and any constants $\kappa\geq0$ and $\lambda>0$, the following function:
$$
\psi(s,x;l_\cdot,\kappa,\lambda)=\exp\Big(\lambda e^{\kappa s}x+\lambda\int_0^sl(r)e^{\kappa r}\dif r\Big),~~~(s,x)\in[0,T]\times\R_+.
$$
\indent (i) In light of (EX), by applying $\rm{It\hat{o}}$-Tanaka's formula to $\psi(t,Y^+_t;\underline{\alpha}_\cdot,\beta,\overline{p})$ and $\psi(t,Y^-_t;\overline{\alpha}_\cdot,\beta,p)$ and using an identical argument as in Proposition 1 of \cite{Fan-Hu-Tang 2020}, we can deduce that for any adapted solution $(Y_\cdot,Z_\cdot)$ of BSDE (\ref{BSDE1.1}), it holds that $\ps$, for each $t\in\T$,
\begin{align}\label{eq:2.3}
     \overline{p} Y^+_t\leq\exp(\overline{p} Y^+_t)\leq\psi(t,Y^+_t;\underline{\alpha}_\cdot,\beta,\overline{p})
     \leq\E\Big[\psi(T,\xi^+;\underline{\alpha}_\cdot,\beta,\overline{p})\big|\F_t\Big],\\
     p Y^-_t\leq\exp(p Y^-_t)\leq\psi(t,Y^-_t;\overline{\alpha}_\cdot,\beta,p)
     \leq\E\Big[\psi(T,\xi^-;\overline{\alpha}_\cdot,\beta,p)\big|\F_t\Big].\nonumber
\end{align}
Then, by virtue of the last two inequalities along with the assumptions of Proposition \ref{cun1}, proceeding as the proof of Proposition 3 in \cite{Briand-Hu 2008} with a localization argument as well as Doob's maximal inequality for martingales, we can conclude that BSDE {\rm{($\ref{BSDE1.1}$)}} admits a solution $(Y_\cdot,Z_\cdot)$ satisfying (\ref{p1}), (\ref{p2}) and $Z_\cdot\in\M^2$. The detailed proof is omitted here.\par
(ii) Since the assumptions in (ii) are stronger than those in (i), it suffices to prove that
$Y_\cdot^+\leq K$ for some constant $K>0$. Indeed, since $\xi^+ +\int_0^T \underline{\alpha}_t \dif t\leq M$ for a constant $M>0$, it follows from \eqref{eq:2.3} that
$$
2\gamma Y^+_t\leq \E[\psi(T,\xi^+;\underline{\alpha}_\cdot,\beta,2\gamma)|\F_t]\leq \exp\left(2\gamma e^{\beta T}M\right),\ \ t\in \T,
$$
which is the desired result. The proof is then complete.\vspace{0.2cm}\hfill\framebox

\begin{remark}
{\rm{Note that under the condition of $\underline{\alpha}_\cdot=\overline{\alpha}_\cdot$, the combination of (i) and (ii) in assumption (EX) can be read as: $\as$, for each $(y,z)\in\R\times\R^d$, we have
$$
-\sgn (y)g(t,y,z)\leq \overline{\alpha}_t+\beta |y|+\frac{\gamma}{2}|z|^2,
$$
which is a form commonly used in the existing literature, for instance, \cite{Fan 2016} and \cite{Fan-Hu-Tang 2020}.}}
\end{remark}

\section{The uniqueness for the case of $g$ satisfying monotonicity in $y$ and convexity in $z$}
\label{section 3}
\setcounter{equation}{0}

This section is devoted to the uniqueness for the unbounded solution of BSDE {\rm{($\ref{BSDE1.1}$)}} with generator $g$ satisfying a monotonicity condition with a general growth in the state variable $y$ and a convexity condition with a quadratic growth in the state variable $z$.

\subsection{Statement of the main result}

We use the following assumptions.\vspace{0.2cm}

\noindent\textbf{(UN1)} $\as$, for each $(y_1,y_2,y,z)\in\R\times\R\times\R\times\R^d$, we have
\begin{enumerate}
\renewcommand{\theenumi}{(1)}
\renewcommand{\labelenumi}{\theenumi}
\item $-\sgn(y_1-y_2)\big(g(t,y_1,z)-g(t,y_2,z)\big)\leq\beta|y_1-y_2|$;
\renewcommand{\theenumi}{(2)}
\renewcommand{\labelenumi}{\theenumi}
\item $g(t,y,\cdot)$ is convex.\vspace{0.1cm}
\end{enumerate}
\noindent\textbf{(UN2)} $\as$, for each $(y,z)\in\R\times\R^d, ~\mathbbm{1}_{y>0}g(t,y,z)\leq \overline{\alpha}_t+\beta|y|+\frac{\gamma}{2}|z|^2$.\vspace{0.2cm}\\
\noindent\textbf{(UN2')} $\as$, for each $(y,z)\in\R\times\R^d$, we have $$\mathbbm{1}_{y>0}\big(g(t,y,z)-g(t,0,z)\big)\leq\overline{\alpha}_t+\beta|y|$$ and
$$|g(t,y,z)-g(t,0,z)|\leq\overline{\alpha}_t+\varphi(|y|).\vspace{0.2cm}$$

If $g(t,y,\cdot)$ is a convex function, then the Legendre-Fenchel transform of $g$ can be defined as follows:
\begin{align}\label{fdy3.1}
f(\omega,t,y,q):=\sup_{z\in\R^d}\big(qz-g(\omega,t,y,z)\big), ~~ (\omega,t,y,q)\in\Omega\times[0,T]\times\R\times\R^{1\times d}.
\end{align}

The following theorem presents several uniqueness results for unbounded solutions of BSDE (\ref{BSDE1.1}), which is the main result of this section.

\begin{theorem}\label{unth3.1}
Suppose that $\xi$ is a terminal value, the generator $g$ is continuous in $(y,z)$ and satisfies {\rm{(EX)}} and {\rm{(UN1)}} , and $\E[\exp(pe^{\beta T}(\xi^-+\int_0^T\overline{\alpha}_t\dif t))]<+\infty$ for some $p>\gamma$.
\begin{enumerate}[leftmargin=0.8cm]
\renewcommand{\theenumi}{{\rm{(i)}}}
\renewcommand{\labelenumi}{\theenumi}
\item If $g$ also satisfies {\rm{(UN2)}}, and $\E[\exp(\overline{p}e^{\beta T}( \xi^++\int_0^T\underline{\alpha}_t\dif t))]<+\infty$ for any $\overline{p}>\gamma$, then BSDE {\rm{($\ref{BSDE1.1}$)}} admits a unique solution $(Y_\cdot,Z_\cdot)$ satisfying {\rm{(\ref{p1})}} and {\rm{(\ref{p2})}} for any $\overline{p}>\gamma$.
\renewcommand{\theenumi}{{\rm{(ii)}}}
\renewcommand{\labelenumi}{\theenumi}
\item If it holds that $\xi^+ +\int_0^T \underline{\alpha}_t \dif t\leq M$ for a constant $M\geq0$, then BSDE {\rm{($\ref{BSDE1.1}$)}} admits a unique solution $(Y_\cdot,Z_\cdot)$ satisfying {\rm{(\ref{p1})}} and $Y_\cdot \leq K$ for a constant $K>0$.
\renewcommand{\theenumi}{{\rm{(iii)}}}
\renewcommand{\labelenumi}{\theenumi}
\item If $g$ also satisfies {\rm{(UN2')}}, and $\E[\exp(\overline{p}e^{\beta T}( \xi^++\int_0^T\underline{\alpha}_t\dif t))]<+\infty$ for any $\overline{p}>\gamma$, then BSDE {\rm{($\ref{BSDE1.1}$)}} admits a unique solution $(Y_\cdot,Z_\cdot)$ satisfying {\rm{(\ref{p1})}} and {\rm{(\ref{p2})}} for any $\overline{p}>\gamma$. Moreover, we have $Y_\cdot=\essinf_{q\in\A}Y_\cdot^q$, and there exists $q^*\in\A$ such that $\as, ~Y_\cdot=Y_\cdot^{q^*}$, where the admissible control set $\A$ is defined as follows:
\begin{align*}
\A:=\bigg\{&(q_s)_{s\in[0,T]}~\text{is an } (\F_t)\text{-progressively measurable~} \R^{1\times d} \text{-valued process}:\\
&\ \ps,\ \int_{0}^{T}|q_s|^2\dif s<+\infty, \ \ ~\E^{\Q}\Big[\int_{0}^{T}|q_s|^2\dif s\Big]<+\infty,\\
&\ \E^{\Q}\Big[|\xi|+\int_{0}^{T}|f(s,0,q_s)|\dif s\Big]<+\infty, \\
&~\text{with}~M_t:=\exp\Big(\int_{0}^{t}q_s\dif B_s-\frac{1}{2}\int_{0}^{t}|q_s|^2\dif s\Big),\ t\in[0,T],\\
&~\text{being a uniformly integrable martingale, and } \frac{\dif \Q}{\dif \mathbb{P}}:=M_{T}\bigg\},\vspace{0.2cm}
\end{align*}
and for each $q\in \A$, $(Y_\cdot^q, Z_\cdot^q)$ is the unique adapted solution of the following BSDE {\rm{(\ref{bsde2})}} such that $Y_\cdot^q$ belongs to class (D) under $\Q$:
\begin{align}\label{bsde2}
Y_t^{q}=\xi+\int_{t}^{T}f(s,Y_s^q,q_s)\dif s+\int_{t}^{T}Z_s^{q}\cdot\dif B_s^q,~~~~t\in[0,T],
\end{align}
with the process $B_t^q:=B_t-\int_0^t q_s\dif s, ~t\in[0,T]$ being a standard $d$-dimensional Brownian motion under $\Q$ and the generator $f$ being defined in {\rm{(\ref{fdy3.1})}}.
\end{enumerate}
\end{theorem}

\begin{remark}\label{rmk-hfsj}
{\rm{With respect to Theorem \ref{unth3.1}, we make the following three comments.
\begin{itemize}
\item[(i)] It is obvious that (UN2) is weaker than (UN2') provided that (EX) holds. The existence of the solution for BSDE (\ref{bsde2}) may not be always true under assumptions (EX), (UN1) and (UN2), which means that we can not adopt the verification argument given in the proof of (iii) of Theorem \ref{unth3.1} to show the uniqueness in (i) and (ii) of Theorem \ref{unth3.1}.

\item[(ii)] In the proof of (i) and (iii) of Theorem \ref{unth3.1}, by further applying the dividing time interval method used in for example the proof of Lemma 3.4 in \cite{Delbaen 2011}, we can verify that if the $\overline{p}$ appearing in (i) and (iii) of Theorem \ref{unth3.1} is replaced with a certain constant $\widetilde{p}>\gamma$, then BSDE (\ref{BSDE1.1}) admits a unique solution $(Y_\cdot,Z_\cdot)$ satisfying (\ref{p1}) and (\ref{p2}) with $\widetilde{p}$ instead of $\overline{p}$, and $Y_\cdot=\essinf_{q\in\overline{\A}}\overline{Y}^{q}_\cdot$ for the case of (iii) of Theorem \ref{unth3.1}, where $\overline{Y}^{q}_\cdot$ and $\overline{\A}$ can be similarly defined as in \cite{Delbaen 2011}.

\item[(iii)] If (1) in assumption (EX) is replaced with
\begin{align*}
-\mathbbm{1}_{y>0}g(t,y,z)\leq \underline{\alpha}_t+\beta|y|+\frac{\gamma}{2}|z|^\alpha, ~\text{where ~} 0<\alpha<2,
\end{align*}
then we can relax the integrability requirement of $\xi^++\int_0^T\underline{\alpha}_t\dif t$ in (i) of Proposition \ref{cun1} along with (i) and (iii) of Theorem \ref{unth3.1} to an $\varepsilon_0 e^{\beta T}$-order exponential moment for some constant $\varepsilon_0>0$ such that BSDE (\ref{BSDE1.1}) admits a (unique) solution $(Y_\cdot,Z_\cdot)$ satisfying (\ref{p1}) and  (\ref{p2}) with $\varepsilon_0$ instead of $\overline{p}$, and $Y_\cdot=\essinf_{q\in\overline{\A}}\overline{Y}^q_\cdot$ for the case of (iii) of Theorem \ref{unth3.1}. In fact, by Young's inequality we know that for each $\varepsilon>0$,
\begin{align*}
\frac{\gamma}{2}|z|^\alpha&=\frac{1}{2}\Big(\frac{2\varepsilon}{\alpha}|z|^2
\Big)^{\frac{\alpha}{2}}\cdot\Big[\gamma^{\frac{2}{2-\alpha}}\Big(\frac{\alpha}
{2\varepsilon}\Big)^{\frac{\alpha}{2-\alpha}}\Big]^{\frac{2-\alpha}{2}}\\
&\leq \frac{\varepsilon}{2}|z|^2+\frac{2-\alpha}{4}\gamma^{\frac{2}{2-\alpha}}
\Big(\frac{\alpha}{2\varepsilon}\Big)^{\frac{\alpha}{2-\alpha}}, ~~\forall z\in\R^d.
\end{align*}
Thus, the desired existence assertion follows by picking $\varepsilon<\overline{p}<\varepsilon_0$ in the proof of Proposition \ref{cun1}, and the uniqueness assertion is also true since we can also fit the proof of Lemma 3.4 in \cite{Delbaen 2011} to our situation, by introducing a positive integer $N$ to divide the time interval $\T$. In particular, we would like to emphasize that these results mentioned above strengthen Corollary 2.2 and Theorem 3.3 in \cite{Delbaen 2011}, where $g$ is Lipschitz continuous in $y$.
\end{itemize}}}
\end{remark}

\begin{example}\label{lz}
{\rm{We give three examples to which Theorem \ref{unth3.1} but no existing result applies.

\begin{enumerate}[leftmargin=1.2cm]
\renewcommand{\theenumi}{{\rm{(i)}}}
\renewcommand{\labelenumi}{\theenumi}
\item For each $(\omega,t,y,z)\in\Omega\times[0,T]\times\R\times\R^d$, define $$g(\omega,t,y,z):=\mathbbm{1}_{y\geq0}\sqrt{y}-\mathbbm{1}_{y<0}y^2
    +\frac{|\sin t|}{2}|z|^2.$$
    It is easy to check that $g$ satisfies neither the locally Lipschitz continuity nor the linear growth with respect to $y$, but it satisfies (EX), (UN1) and (UN2) with $\overline{\alpha}_t=1, ~\underline{\alpha}_t=1, ~\beta=1, ~\gamma=1$ and $\varphi(u)=u^2+\sqrt{u}$. It then follows from (i) of Theorem \ref{unth3.1} that for each $\xi$ such that $\xi^-$ admits an exponential moment of order $pe^T$ with $p>1$ and $\xi^+$ admits all exponential moments, BSDE $(\xi,T,g)$ admits a unique solution $(Y_\cdot,Z_\cdot)$ satisfying (\ref{p1}) and (\ref{p2}) for any $\overline{p}>1$.

\renewcommand{\theenumi}{{\rm{(ii)}}}
\renewcommand{\labelenumi}{\theenumi}
\item For each $(\omega,t,y,z)\in\Omega\times[0,T]\times\R\times\R^d$, define $$g(\omega,t,y,z):=\mathbbm{1}_{y\geq0}e^y+\mathbbm{1}_{y<0}
    (\sqrt[3]{y}+y^3+1)+\frac{1}{2}(|z|^{\frac{3}{2}}\vee |z|^2).$$
    Clearly, $g$ satisfies neither the locally Lipschitz continuity nor the linear growth with respect to $y$, but it satisfies (EX) and (UN1) with $\overline{\alpha}_t=2, ~\underline{\alpha}_t=2, ~\beta=1, ~\gamma=1$ and $\varphi(u)=e^u+\sqrt[3]{u}+u^3$. It then follows from (ii) of Theorem \ref{unth3.1} that for each $\xi$ such that $\xi^-$ admits an exponential moment of order $pe^T$ with $p>1$ and $\xi^+$ is bounded, BSDE $(\xi,T,g)$ admits a unique solution $(Y_\cdot,Z_\cdot)$ satisfying (\ref{p1}) and $Y_\cdot\leq K$ for some constant $K>0$.

\renewcommand{\theenumi}{{\rm{(iii)}}}
\renewcommand{\labelenumi}{\theenumi}
\item For each $(\omega,t,y,z)\in\Omega\times[0,T]\times\R\times\R^d$, define $$g(\omega,t,y,z):=\mathbbm{1}_{y\geq1}\sqrt{y-1}+\mathbbm{1}_{y<1}(y-1)^3
    +\frac{1}{2}|z|^2.$$
    This is just the generator defined in \eqref{eq:1.4}. It is easy to check that $g$ satisfies neither the locally Lipschitz continuity nor the linear growth with respect to $y$, but it satisfies (EX), (UN1) and (UN2') with $\overline{\alpha}_t=1, ~\underline{\alpha}_t=1, ~\beta=1, ~\gamma=1$ and $\varphi(u)=1+u^3$. It then follows from (iii) of Theorem \ref{unth3.1} that for each $\xi$ such that $\xi^-$ admits an exponential moment of order $pe^{T}$ with $p>1$ and $\xi^+$ admits all exponential moments, BSDE $(\xi,T,g)$ admits a unique solution $(Y_\cdot,Z_\cdot)$ satisfying (\ref{p1}) and (\ref{p2}) for any $\overline{p}>1$. Moreover, we have $Y_\cdot=\essinf_{q\in\A}Y_\cdot^q$, where $(Y_\cdot^q, Z_\cdot^q)$ is the unique solution of BSDE {\rm{(\ref{bsde2})}} such that $Y_\cdot^q$ belongs to class (D) under $\Q$.
    \end{enumerate}
    }}
\end{example}

At the end of this subsection, it should be especially emphasized that by the analysis presented in the introduction we can conclude that for the generator $g$ defined in (iii) of Example \ref{lz}, neither $f(t,y,z):=-g(t,y,-z)$ nor $f(t,y,z):=g(t,-y,z)$ satisfies the extended convexity condition \eqref{eq:1.3} proposed in \cite{Fan-Hu-Tang 2020}. Furthermore, in a similar manner it can be proved that this assertion also holds for the generator $g$ defined in (i) and (ii) of Example \ref{lz}, respectively. As a result, those conclusions stated in Example \ref{lz} can not be obtained by the results proposed in \cite{Fan-Hu-Tang 2020}.

\subsection{Proof of the main result}

The existence has been proved in Proposition \ref{cun1}. Now let us verify the uniqueness.\vspace{0.2cm}

\noindent\textbf{Proof of (i) of Theorem \ref{unth3.1}.} According to (2) of (EX), (UN2) and (UN1), we can deduce that the generator $f$ defined in (\ref{fdy3.1}) satisfies the following properties: $\as$,
\begin{itemize}
\item~~$\forall (y,q)\in\R\times{\R^{1\times d}}, ~f(t,y,q)\geq-\overline{\alpha}_t-\beta |y|+\frac{1}{2\gamma}|q|^2;$
\item~~$\forall q\in{\R^{1\times d}}, ~\forall(y_1,y_2)\in{\R^2},~\sgn(y_1-y_2)\big(f(t,y_1,q)-f(t,y_2,q)\big)\leq \beta|y_1-y_2|$;
\item~~$\forall y\in\R, ~f(t,y,\cdot)$ is convex.
\end{itemize}
The last two assertions are obvious. We only need to verify the first assertion. In fact, it follows from (2) of (EX) and (UN2) that $\as, ~\forall (y,z)\in\R\times\R^d$,
$$
g(t,y,z)\leq \overline{\alpha}_t+\beta|y|+\frac{\gamma}{2}|z|^2.
$$
Then, by (\ref{fdy3.1}) we have $\as$, for each $(y,q)\in\R\times\R^{1\times d}$,
$$
f(t,y,q)\geq -\overline{\alpha}_t-\beta|y|+\sup_{z\in\R^d}\big(qz-\frac{\gamma}{2}|z|^2\big)
=-\overline{\alpha}_t-\beta|y|+\frac{1}{2\gamma}|q|^2.
$$

Let us proceed the proof by giving a uniform integrability result.

\begin{proposition}\label{pro-Q}
Assume that the generator $g$ satisfies assumptions {\rm{(EX)}}, {\rm{(UN1)}} and {\rm{(UN2)}}. Let $(Y_\cdot, Z_\cdot)$ be an adapted solution of BSDE {\rm{($\ref{BSDE1.1}$)}} such that $Y_\cdot$ satisfies {\rm{(\ref{p1})}} and {\rm{(\ref{p2})}} for each $\overline{p}>\gamma$. Then, for all $(\F_t)$-progressively measurable processes $(q^*_s)_{s\in[0,T]}$ valued in $\R^{1\times d}$ such that $q_s^*\in\partial_z g(s,Y_s,Z_s)$ for all $s\in[0,T]$, $\mathcal{E}(q^*)$ is a uniformly integrable martingale and defines a probability $\Q^*\sim\mathbb{P}$. Moreover, we have
\begin{align}\label{h}
\E^{\Q^*}\Big[\int_0^T|q^*_s|^2\dif s\Big]<+\infty.
\end{align}
\end{proposition}

\noindent\textbf{Proof.}
In view of $q_s^*\in\partial_z g(s,Y_s,Z_s)$, we have $g(s,Y_s,Z_s)=q_s^*Z_s-f(s,Y_s,q_s^*)$, which along with $f(t,y,q)\geq-\overline{\alpha}_t-\beta |y|+\frac{1}{2\gamma}|q|^2$ yields
\begin{align*}
g(s,Y_s,Z_s)=q_s^*Z_s-f(s,Y_s,q_s^*)&\leq \frac{1}{2}\Big(\frac{1}{2\gamma}|q_s^*|^2+2\gamma|Z_s|^2\Big)+\overline{\alpha}_s
+\beta|Y_s|-\frac{1}{2\gamma}|q_s^*|^2\\
&=\gamma|Z_s|^2+\overline{\alpha}_s+\beta|Y_s|-\frac{1}{4\gamma}|q_s^*|^2;
\end{align*}
\vspace{-0.6cm}
\begin{align*}
\frac{1}{4\gamma}|q_s^*|^2&\leq-g(s,Y_s,Z_s)+\gamma|Z_s|^2+\overline{\alpha}_s+\beta|Y_s|.
\end{align*}
It follows from the last inequality that $\ps,~\int_0^T|q_s^*|^2\dif s<+\infty$.\par
For each $n\geq1$, let us define
$$
M_t:=\exp\Big(\int_{0}^{t}q_s^*\dif B_s-\frac{1}{2}\int_{0}^{t}|q_s^*|^2\dif s\Big), ~t\in\T,
$$
$$
\tau_n:=\inf \Big\{t \in[0, T]: \int_0^t|q_s^*|^2 \dif s+\int_0^t|Z_s|^2 \dif s\geq n\Big\} \wedge T, \quad \frac{\dif \Q_n^*}{\dif \mathbb{P}}:=M_{\tau_n}.
$$
We set $B_t^{q^*}:=B_t-\int_0^t q^*_s\dif s, ~t\in\T$. Then $(B_t^{q^*})_{t\in [0,\tau_n]}$ is a standard $d$-dimensional Brownian motion under $\Q_n^*$. Now we show that $(M_{\tau_n})_n$ is uniformly integrable.

According to (\ref{p1}) and the Fenchel's inequality, we have
\begin{align*}
\E^{\Q_n^{*}}\Big[\sup_{t\in\T} A_\beta(t)\Big]&=\E\Big[M_{\tau_n} \sup_{t\in\T} A_\beta(t)\Big] \\
&\leq \E\Big[\exp\big(p\sup_{t\in\T} A_\beta(t)\big)\Big]+\frac{1}{p}\E\big[M_{\tau_n}(\ln M_{\tau_n}-\ln p-1)\big]\\
&=\E\Big[\exp\big(p\sup_{t\in\T} A_\beta(t)\big)\Big]+\frac{1}{p}\E[M_{\tau_n}\ln M_{\tau_n}]-\frac{1}{p}(\ln p+1).
\end{align*}
Some uncomplicated calculations give
\begin{align}\label{t}
\notag\E[M_{\tau_n}\ln {M_{\tau_n}}]&=\E^{\Q_n^{*}}[\ln M_{\tau_n}]=\E^{\Q_n^{*}}\Big[\int_{0}^{\tau_n}q^*_s\dif B_s-\frac{1}{2}\int_{0}^{\tau_n}|q^*_s|^2\dif s\Big]\\
&=\E^{\Q_n^{*}}\Big[\int_{0}^{\tau_n}q^*_s\dif B^{q^*}_s+\frac{1}{2}\int_{0}^{\tau_n}|q^*_s|^2\dif s\Big]=\frac{1}{2}\E^{\Q_n^{*}}\Big[\int_{0}^{\tau_n}|q^*_s|^2\dif s\Big].
\end{align}
The calculation similar to (\ref{t}) will be used several times later, and this calculation process will not be repeated in detail. For ease of notations, we denote the sum of constants by $C_p$, then we have
\begin{align}\label{3a}
\E^{\Q_n^{*}}\Big[\sup_{t\in\T} A_\beta(t)\Big]\leq C_p+\frac{1}{2p} \E^{\Q_n^*}\Big[\int_0^{\tau_n}|q_s^*|^2 \dif s\Big],
\end{align}
and for each $\overline{p}>\gamma$, in the same manner,
\begin{align}\label{3aa}
\E^{\Q_n^{*}}\Big[\sup_{t\in\T} \overline{A}_\beta(t)\Big]\leq C_{\overline{p}}+\frac{1}{2\overline{p}} \E^{\Q_n^*}\Big[\int_0^{\tau_n}|q_s^*|^2 \dif s\Big].
\end{align}

Since $g(s,Y_s,Z_s)=q_s^*Z_s-f(s,Y_s,q_s^*)$, by BSDE (\ref{BSDE1.1}) we have
$$
Y_{t\wedge\tau_n}=Y_{\tau_n}+\int_{t\wedge\tau_n}^{\tau_n}f(s,Y_s,q_s^*)\dif s+\int_{t\wedge\tau_n}^{\tau_n} Z_s \cdot \dif B_s^{q^*}, ~~t\in\T.
$$
Applying $\rm{It\hat{o}}$'s formula to $e^{\beta s} Y_s$ gives
\begin{align*}
\dif (e^{\beta s} Y_s)=e^{\beta s}\big(\beta Y_s-f(s,Y_s,q_s^*)\big)\dif s-e^{\beta s}Z_s\cdot\dif B_s^{q^*}, ~s\in[0,\tau_n].
\end{align*}
In view of $f(t,y,q)\geq-\overline{\alpha}_t-\beta|y|+\frac{1}{2\gamma}|q|^2$ and $Y_{\tau_n}\geq-Y_{\tau_n}^-$, we deduce that
\begin{align*}
Y_0&=e^{\beta \tau_n} Y_{\tau_n}+\int_0^{\tau_n}e^{\beta s}\big(f(s,Y_s,q_s^*)-\beta Y_s\big)\dif s+\int_0^{\tau_n}e^{\beta s}Z_s\cdot\dif B_s^{q^*}\\
&\geq -e^{\beta \tau_n} Y_{\tau_n}^-+\int_0^{\tau_n}e^{\beta s}\Big(-\overline{\alpha}_s-\beta|Y_s|+\frac{1}{2\gamma}|q_s^*|^2-\beta Y_s\Big)\dif s+\int_0^{\tau_n}e^{\beta s}Z_s\cdot\dif B_s^{q^*}\\
&= -e^{\beta \tau_n} Y_{\tau_n}^-+\int_0^{\tau_n}e^{\beta s}\Big(-\overline{\alpha}_s-2\beta Y_s^++\frac{1}{2\gamma}|q_s^*|^2\Big)\dif s+\int_0^{\tau_n}e^{\beta s}Z_s\cdot\dif B_s^{q^*}.
\end{align*}
It then follows from (\ref{3a}) and (\ref{3aa}) that for each $\overline{p}>\gamma$,
\begin{align}\label{tt}
\notag Y_0&\geq \E^{\Q_n^*}\Big[-e^{\beta \tau_n} Y_{\tau_n}^-+\int_0^{\tau_n}e^{\beta s}\Big(-\overline{\alpha}_s-2\beta Y_s^++\frac{1}{2\gamma}|q_s^*|^2\Big)\dif s\Big]\\
\notag &\geq -\E^{\Q_n^*}\Big[e^{\beta \tau_n}\Big( Y_{\tau_n}^-+\int_0^{\tau_n}\overline{\alpha}_s\dif s\Big)\Big]-2\beta \E^{\Q_n^*}\Big[\int_0^{\tau_n}e^{\beta s}Y_s^+\dif s\Big]+\frac{1}{2\gamma}\E^{\Q_n^*}\Big[\int_0^{\tau_n}e^{\beta s}|q_s^*|^2\dif s\Big]\\
\notag &\geq-\E^{\Q_n^*}\Big[A_\beta(\tau_n)\Big]-2\beta T\E^{\Q_n^*}\Big[\sup_{s\in\T}e^{\beta s}Y_s^+\Big]+\frac{1}{2\gamma}\E^{\Q_n^*}\Big[\int_0^{\tau_n}|q_s^*|^2\dif s\Big]\\
&\geq-C_p-2\beta T C_{\overline{p}}+\frac{1}{2}\Big(\frac{1}{\gamma}-\frac{1}{p}-\frac{2\beta T}{\overline{p}}\Big)\E^{\Q_n^*}\Big[\int_0^{\tau_n}|q_s^*|^2\dif s\Big].
\end{align}
Set $\overline{p}>\frac{2\beta Tp\gamma}{p-\gamma}$. Then $\frac{1}{2}\Big(\frac{1}{\gamma}-\frac{1}{p}-\frac{2\beta T}{\overline{p}}\Big)>0.$ Hence, in light of (\ref{t}) and (\ref{tt}), we have
\begin{align}\label{3g}
2\E[M_{\tau_n}\ln M_{\tau_n}]=\E^{\Q_n^{*}}\Big[\int_0^{\tau_n}|q_s^*|^2\dif s\Big]<C_{p,\overline{p},T,\beta,\gamma},
\end{align}
where $C_{p,\overline{p},T,\beta,\gamma}$ is a positive constant independent of $n$. Hence,  $\sup_n\E[M_{\tau_n}\ln M_{\tau_n}]<+\infty$. Thus, $(M_{\tau_n})_n$ is uniformly integrable by the de La Vall$\acute{\rm{e}}$e Poussin lemma, and then $\E[M_T]=1$, which means that $\mathcal{E}(q^*)=(M_t)_{t\in\T}$ is a uniformly integrable martingale and defines a probability $\Q^*\sim\mathbb{P}$. Moreover, applying Fatou's lemma in (\ref{3g}), we obtain
$$
2\E[M_T\ln M_T]=\E^{\Q^*}\Big[\int_0^T|q_s^*|^2\dif s\Big]\leq\liminf_n\E^{\Q_n^* }\Big[\int_0^{\tau_n}|q_s^*|^2\dif s\Big]<+\infty.
$$
The proof of Proposition \ref{pro-Q} is then complete.\hfill\framebox

\vspace{0.2cm}

\indent In the sequel, let $(Y_\cdot, Z_\cdot)$ and $(Y'_\cdot,Z'_\cdot)$ be two adapted solutions of BSDE {\rm{($\ref{BSDE1.1}$)}} such that for each $\overline{p}>\gamma$,
\begin{align}\label{h1}
\E\Big[\exp\big(p\sup_{t\in[0,T]}A_\beta(t)\big)\Big]<+\infty, ~ \E\Big[\exp\big(\overline{p}\sup_{t\in[0,T]}\overline{A}_\beta(t)\big)\Big]<+\infty,
\end{align}
\vspace{-0.5cm}
\begin{align}\label{h2}
\E\Big[\exp\big(p\sup_{t\in[0,T]}A'_\beta(t)\big)\Big]<+\infty, ~ \E\Big[\exp\big(\overline{p}\sup_{t\in[0,T]}\overline{A'}_\beta(t)\big)\Big]<+\infty,
\end{align}
with $A'_\beta(t):=e^{\beta t}\big((Y'_t)^-+\int_0^t\overline{\alpha}_s\dif s\big)$ and $\overline{A'}_\beta(t):=e^{\beta t}\big((Y'_t)^++\int_0^t\underline{\alpha}_s\dif s\big)$ as well as $A_\beta(\cdot)$ and $\overline{A}_\beta(\cdot)$ being defined respectively in (\ref{p1}) and (\ref{p2}).

In order to obtain the desired uniqueness, by a symmetry argument it suffices to show that for each $t\in[0,T], ~\ps,~Y_t\geq Y'_t$. For $t\in[0,T)$, let us denote $A:=\{Y_t<Y'_t\}$ and the stopping time $\tau:=\inf\{s\geq t | Y_s\geq Y'_s\}$. Then for $s\in[t,\tau]$, we have $\mathbbm{1}_A Y_s\leq \mathbbm{1}_A Y'_s$ and $\mathbbm{1}_A Y_\tau=\mathbbm{1}_A Y'_\tau$ since $t\rightarrow Y_t$ is continuous. It implies that $(\mathbbm{1}_A(Y_s-Y'_s))_{s\in[t,\tau]}$ is a non-positive process.

Let us consider an $(\F_s)$-progressively measurable process $(q_s^*)_{s\in[0,T]}$ valued in $\R^{1\times d}$ such that $q_s^*\in\partial_z g(s,Y_s,Z_s)$ for all $s\in[0,T]$. Thanks to Proposition \ref{pro-Q}, we can conclude that $(M_t^*)_{t\in\T}$ defines a probability measure $\Q^*\thicksim \mathbb{P}$, where
$$
M_t^*:=\exp\Big(\int_{0}^{t}q_s^*\dif B_s-\frac{1}{2}\int_{0}^{t}|q_s^*|^2\dif s\Big)~~\text{and}~~\frac{\dif \Q^*}{\dif \mathbb{P}}:=M_T^*
$$
We set $B_t^{q^*}:=B_t-\int_0^t q^*_s\dif s,~t\in\T$. Thanks to Girsanov's theorem, $(B_t^{q^*})_{t\in\T}$ is a standard $d$-dimensional Brownian motion under $\Q^*$. Then,
$$
\dif \big(Y_s-Y'_s\big)=\Big(g(s,Y_s,Z_s)-g(s,Y'_s,Z'_s)-q_s^*(Z_s-Z'_s)\Big)\dif s-(Z_s-Z'_s)\cdot\dif B_s^{q^*}, ~s\in\T.
$$
Applying $\rm{It\hat{o}}$'s formula to $e^{\beta s}\mathbbm{1}_A(Y_s-Y'_s)$, we obtain that for each $s\in[t,\tau]$,
\begin{align*}
\dif \big[e^{\beta s}\mathbbm{1}_A(Y_s-Y'_s)\big]=&\beta e^{\beta s}\mathbbm{1}_A(Y_s-Y'_s)\dif s+e^{\beta s}\mathbbm{1}_A\big(g(s,Y_s,Z_s)-g(s,Y_s,Z'_s)-q_s^*(Z_s-Z'_s)\\
&+g(s,Y_s,Z'_s)-g(s,Y'_s,Z'_s)\big)\dif s-e^{\beta s}\mathbbm{1}_A(Z_s-Z'_s)\cdot\dif B_s^{q^*}.
\end{align*}
Noticing by $q_s^*\in\partial_z g(s,Y_s,Z_s)$ that $g(s,Y_s,Z_s)-g(s,Y_s,Z'_s)-q_s^*(Z_s-Z'_s)\leq0$, we have
\begin{align*}
&\dif \big[e^{\beta s}\mathbbm{1}_A(Y_s-Y'_s)\big]\\
&\ \ \ \leq e^{\beta s}\mathbbm{1}_A\big(\beta(Y_s-Y'_s)+g(s,Y_s,Z'_s)-g(s,Y'_s,Z'_s)\big)\dif s-e^{\beta s}\mathbbm{1}_A(Z_s-Z'_s)\cdot\dif B_s^{q^*}.
\end{align*}
In light of (1) of (UN1) along with the fact that $(\mathbbm{1}_A(Y_s-Y'_s))_{s\in[t,\tau]}$ is a non-positive process, it is not very hard to verify that
$$
\mathbbm{1}_A\big(\beta(Y_s-Y'_s)+g(s,Y_s,Z'_s)-g(s,Y'_s,Z'_s)\big)\leq0, ~~s\in[t,\tau].
$$
Then, we have
$$
\dif \big[e^{\beta s}\mathbbm{1}_A(Y_s-Y'_s)\big]\leq-e^{\beta s}\mathbbm{1}_A(Z_s-Z'_s)\cdot\dif B_s^{q^*}, ~~s\in[t,\tau].
$$
Define the following stopping time, for each $n\geq1$,
$$
\sigma^t_n:=\inf\Big\{s\geq t:\int_t^s|Z_u|^2\dif u+\int_t^s|Z'_u|^2\dif u\geq n\Big\}\wedge\tau
$$
with the convention $\inf\emptyset=+\infty$. Then, we obtain
$$
e^{\beta t}\mathbbm{1}_A(Y_t-Y'_t)\geq e^{\beta \sigma^t_n}\mathbbm{1}_A(Y_{\sigma^t_n}-Y'_{\sigma^t_n})+\int_t^{\sigma^t_n}e^{\beta s}\mathbbm{1}_A(Z_s-Z'_s)\cdot\dif B_s^{q^*}.
$$
It follows that
\begin{align}\label{j}
e^{\beta t}\mathbbm{1}_A(Y_t-Y'_t)\geq \E^{\Q^*}\Big[e^{\beta \sigma^t_n}\mathbbm{1}_A(Y_{\sigma^t_n}-Y'_{\sigma^t_n})\big|\F_t\Big].
\end{align}

For ease of notations, let us set $\widehat{Y}_\cdot:=Y_\cdot-Y'_\cdot$, so $|Y_{\sigma^t_n}-Y'_{\sigma^t_n}|\leq\sup_{s\in\T}|\widehat{Y}_s|$. In light of (\ref{h1}) and (\ref{h2}) along with (\ref{h}) in Proposition \ref{pro-Q}, by using an identical proof to (\ref{3a}) and (\ref{3aa}) we obtain that, picking $\overline{p}=p$,
$$
\E^{\Q^*}\Big[\sup_{s\in\T}Y_s^-\Big]\leq\E^{\Q^*}\Big[\sup_{s\in\T}A_\beta(s)\Big]\leq C_p+\frac{1}{2p} \E^{\Q^*}\Big[\int_0^T|q_s^*|^2 \dif s\Big]=:\overline{C}_p,\vspace{-0.3cm}
$$
\begin{align*}
\E^{\Q^*}\Big[\sup\limits_{s\in\T}Y_s^+\Big]\leq \overline{C}_{p}, ~ \E^{\Q^*}\Big[\sup\limits_{s\in\T}(Y'_s)^-\Big]\leq \overline{C}_p ~~\text{and~~} \E^{\Q^*}\Big[\sup\limits_{s\in\T}(Y'_s)^+\Big]\leq \overline{C}_{p}.
\end{align*}
Then, we have
\begin{align*}
\E^{\Q^*}\Big[\sup_{s\in\T}|\widehat{Y}_s|\Big]&=\E^{\Q^*}\Big[\sup_{s\in\T}
(\widehat{Y}^+_s+\widehat{Y}^-_s)\Big]
\leq 4C_p+\frac{2}{p} \E^{\Q^*}\Big[\int_0^T|q_s^*|^2 \dif s\Big]<+\infty.
\end{align*}
Consequently, by sending $n\rightarrow\infty$ in (\ref{j}) and applying the Lebesgue's dominated convergence theorem we get that
\begin{align*}
\ps,\ \ e^{\beta t}\mathbbm{1}_A(Y_t-Y'_t)\geq \E^{\Q^*}\Big[e^{\beta \tau}(\mathbbm{1}_A Y_\tau-\mathbbm{1}_A Y'_\tau)\big|\F_t\Big]=0,
\end{align*}
which implies that for each $t\in\T, ~\ps, ~\mathbbm{1}_A(Y_t-Y'_t)\geq 0$. Considering that $(\mathbbm{1}_A(Y_s-Y'_s))_{s\in[t,\tau]}$ is a non-positive process, we know that for each $t\in\T, ~\ps,~ \mathbbm{1}_A(Y_t-Y'_t)=0$, that is to say, $\mathbb{P}(A)=0$ and $\ps,~Y_t\geq Y'_t$. The proof of (i) of Theorem \ref{unth3.1} is complete.\hfill\framebox

\vspace{0.3cm}

\noindent\textbf{Proof of (ii) of Theorem \ref{unth3.1}.}
According to (2) and (3) of (EX) and (UN1), we can deduce that the function $f$ defined in (\ref{fdy3.1}) satisfies the following properties: $\as$,
\begin{itemize}
\item~~$\forall (y,q)\in\R_+\times{\R^{1\times d}}, ~f(t,y,q)\geq-\overline{\alpha}_t-\varphi(y^+)+\frac{1}{2\gamma}|q|^2;$
\item~~$\forall (y,q)\in\R_-\times{\R^{1\times d}}, ~f(t,y,q)\geq -\overline{\alpha}_t-\beta y^-+\frac{1}{2\gamma}|q|^2;$
\item~~$\forall q\in{\R^{1\times d}}, ~\forall(y_1,y_2)\in{\R^2},~\sgn(y_1-y_2)\big(f(t,y_1,q)-f(t,y_2,q)\big)\leq \beta|y_1-y_2|$;
\item~~$\forall y\in\R, ~f(t,y,\cdot)$ is convex.
\end{itemize}
It then follows that $\as$, for each $y\leq K$ with some constant $K>0$, we have
\begin{align*}
f(t,y,q)\geq -\overline{\alpha}_t-\beta y^--\varphi(y^+)+\frac{1}{2\gamma}|q|^2
\geq-\overline{\alpha}_t-\varphi(K)-\beta y^-+\frac{1}{2\gamma}|q|^2.
\end{align*}
With the above observation along with the requirement that $Y_\cdot \leq K$, the desired uniqueness assertion can be proved in the same way as in (i) of Theorem \ref{unth3.1}.\hfill\framebox

\vspace{0.3cm}

\noindent\textbf{Proof of (iii) of Theorem \ref{unth3.1}.}
Note by (3) of (EX) that $\as, ~|g(t,0,z)|\leq\overline{\alpha}_t+\frac{\gamma}{2}|z|^2$ for each $z\in\R^d$. Assumptions (EX) and (UN2') can imply (UN2) with $2\overline{\alpha}_\cdot$ instead of $\overline{\alpha}_\cdot$. Therefore, the assumptions in (iii) of Theorem \ref{unth3.1} are stronger than those in (i) of Theorem \ref{unth3.1}. Then, we can conclude that $f$ defined in (\ref{fdy3.1}) satisfies the following properties: $\as$,
\begin{itemize}
\item[(1)] $\forall (y,q)\in\R\times{\R^{1\times d}}, ~f(t,y,q)\geq-\overline{\alpha}_t-\beta|y|+\frac{|q|^2}{2\gamma};$
\item[(2)] $\forall q\in{\R^{1\times d}},~ ~\forall(y_1,y_2)\in{\R^2},~\sgn(y_1-y_2)\big(f(t,y_1,q)-f(t,y_2,q)\big)\leq \beta|y_1-y_2|$;
\item[(3)] $\forall y\in\R, ~f(t,y,\cdot)$ is convex.
\end{itemize}
Furthermore, it follows from (2) that $\as$, for each $(y,q)\in\R_-\times{\R^{1\times d}}, ~f(t,0,q)\leq f(t,y,q)+\beta|y|$. By (UN2') and (\ref{fdy3.1}) we can deduce that $\as$,
\begin{itemize}
\item[(4)] $\forall (y,q)\in\R^+\times{\R^{1\times d}}, ~f(t,0,q)\leq f(t,y,q)+\overline{\alpha}_t+\beta|y|;$
\item[(5)] $\forall (y,q)\in\R\times{\R^{1\times d}}, ~|f(t,y,q)-f(t,0,q)|\leq \overline{\alpha}_t+\varphi(|y|).$
\end{itemize}
Therefore, $\as$, for each $(y,q)\in\R\times{\R^{1\times d}}$, we have
\begin{align}\label{x}
f(t,0,q)\leq f(t,y,q)+\overline{\alpha}_t+\beta|y|.
\end{align}
Now, let $q$ be in $\A$, if this set is not empty. In light of (2), (5), the definition of $\A$ and  Girsanov's theorem, by Theorem 6.2 in \cite{Briand 2003} we know that BSDE {\rm{(\ref{bsde2})}} admits a unique solution $(Y_\cdot^q,Z_\cdot^q)$ such that $(Y_t^q)_{t\in[0,T]}$ belongs to class $(D)$ under $\Q$. Let $(Y_\cdot,Z_\cdot)$ be a solution of BSDE (\ref{BSDE1.1}) satisfying (\ref{p1}) and (\ref{p2}) for any $\overline{p}>\gamma$. We need to prove
$$
Y_\cdot=\essinf\limits_{q\in\A}Y_\cdot^q,
$$
which consists in minimizing $Y_\cdot^{q}$ among all the admissible controls $q\in\A$. The following proof is divided into three steps.\vspace{0.2cm}

\textbf{First step.} Let us start with showing that $Y_\cdot\leq Y_\cdot^{q}$ for any $q\in\A$. Firstly, we have $Y_T=Y^{q}_T=\xi$. For each $t\in[0,T]$ and each integer $m\geq1$, denote the following stopping time
 $$\tau^t_m:=\inf\Big\{s\geq t,\int_t^s|Z_u|^2\dif u+\int_t^s|Z^{q}_u|^2\dif u+\int_t^s|q_u|^2\dif u\geq m\Big\}\wedge T$$
with the convention $\inf\emptyset=+\infty$. Applying $\rm{It\hat{o}}$-Tanaka's formula to $e^{\beta s} (Y_s-Y_s^{q})^+$, we obtain
\begin{align}\label{q}
\dif \big[e^{\beta s}(Y_s-Y_s^q)^+\big]\geq &\Big(\beta e^{\beta s}(Y_s-Y_s^q)^++e^{\beta s}\mathbbm{1}_{\{Y_s-Y_s^q>0\}}\big(f(s,Y_s^q,q_s)+g(s,Y_s,Z_s)-q_sZ_s\big)\Big)\dif s \notag \\
&-e^{\beta s}\mathbbm{1}_{\{Y_s-Y_s^q>0\}}(Z_s-Z^{q}_s)\cdot\dif B_s^q, ~~~s\in[t,\tau].
\end{align}
By (1) of (UN1) and (\ref{fdy3.1}), we have
\begin{align*}
\mathbbm{1}_{\{Y_s-Y_s^q>0\}}&\Big(f(s,Y_s^q,q_s)+g(s,Y_s,Z_s)-q_sZ_s\Big)\\
=&\mathbbm{1}_{\{Y_s-Y_s^q>0\}}\Big(\big(f(s,Y_s^q,q_s)+g(s,Y_s^q,Z_s)-q_sZ_s\big)
+\big(g(s,Y_s,Z_s)-g(s,Y_s^q,Z_s)\big)\Big)\\
\geq &\mathbbm{1}_{\{Y_s-Y_s^q>0\}}\big(0-\beta(Y_s-Y_s^q)^+\big)\geq-\beta(Y_s-Y_s^q)^+, ~~s\in[t,\tau],
\end{align*}
which along with (\ref{q}) yields that
\begin{align*}
\dif \big[e^{\beta s}(Y_s-Y_s^q)^+\big]
\geq -e^{\beta s}\mathbbm{1}_{\{Y_s-Y_s^q>0\}}(Z_s-Z^{q}_s)\cdot\dif B_s^q, ~~s\in[t,\tau].
\end{align*}
Hence, for each $t\in\T$ and $m\geq1$, we have
$$
e^{\beta t} (Y_t-Y_t^{q})^+\leq e^{\beta \tau^t_m} (Y_{\tau^t_m}-Y_{\tau^t_m}^{q})^++\int_t^{\tau^t_m}e^{\beta s}\mathbbm{1}_{\{Y_s-Y_s^q>0\}}(Z_s-Z^{q}_s)\cdot\dif B_s^q.
$$
It then follows that
\begin{align}\label{b}
e^{\beta t} (Y_t-Y_t^{q})^+\leq\E^{\Q}\big[e^{\beta\tau^t_m} (Y_{\tau^t_m}-Y_{\tau^t_m}^{q})^+\big|\F_t\big].
\end{align}
In light of (\ref{p1}) and (\ref{p2}), by picking $\overline{p}=p$, the Fenchel's inequality gives
\begin{align*}
\E^{\Q}\Big[\sup_{s\in\T} Y_s^-+\sup_{s\in\T} Y_s^+\Big]&=\E\Big[M_T \sup_{s\in\T} Y_s^-\Big]+\E\Big[M_T \sup_{s\in\T} Y_s^+\Big]\\
&\leq \E\Big[\exp\Big(p\sup_{s\in\T} Y_s^+\Big)\Big]+\E\Big[\exp\Big(p\sup_{s\in\T} Y_s^-\Big)\Big]\\
&~~~~~+\frac{2}{p}\E[M_T(\ln M_T-\ln p-1)]\\
&\leq 2C_p+\frac{1}{p}\E^{\Q}\Big[\int_{0}^{T}|q_s|^2\dif s\Big]<+\infty.
\end{align*}
Thus, in light of the last inequality and the fact that $Y_\cdot^q$ belongs to class $(D)$ under $\Q$, by letting $m\rightarrow\infty$ in (\ref{b}) we obtain
\begin{align*}
\ps, ~~e^{\beta t} (Y_t-Y_t^{q})^+\leq\E^{\Q}\big[e^{\beta T} (Y_{T}-Y_{T}^{q})^+\big|\F_t\big]=0
\end{align*}
which means that for each $t\in\T, ~\ps, ~Y_t\leq Y_t^q$.
\vspace{0.2cm}

\noindent\textbf{Second step.} Set $q_s^*\in\partial_{z}g(s,Y_s,Z_s)$. Then we have
$$
f(s,Y_s,q_s^*)=q_s^*Z_s-g(s,Y_s,Z_s), ~~s\in[0,T].
$$
Consequently, if $q_\cdot^*\in\A$, then by BSDEs (\ref{BSDE1.1}) and (\ref{bsde2}) along with uniqueness of the integrable solution of BSDE (\ref{bsde2}), we can conclude that $Y_\cdot=Y_\cdot^{q^*}$.\vspace{0.2cm}

\noindent\textbf{Third step.} We conclude the proof by verifying $q^*\in\A$. Note that the assumptions in (iii) of Theorem \ref{unth3.1} are stronger than those in (i) of Theorem \ref{unth3.1}. By Proposition \ref{pro-Q} we know that $\ps,~\int_0^T|q_s^*|^2\dif s<+\infty$, $(M^*_t)_{t\in\T}$ is a uniformly integrable martingale and defines a probability measure $\Q^*\sim\mathbb{P}$, where
$$
M^*_t:=\exp\Big(\int_{0}^{t}q_s^*\dif B_s-\frac{1}{2}\int_{0}^{t}|q_s^*|^2\dif s\Big), ~t\in\T ~~\text{and}~~\frac{\dif \Q^*}{\dif \mathbb{P}}:=M^*_T.
$$
We set $B_t^{q^*}:=B_t-\int_0^t q^*_s\dif s,~t\in\T$. Thanks to Girsanov's theorem, $(B_t^{q^*})_{t\in\T}$ is a standard $d$-dimensional Brownian motion under $\Q^*$. Moreover, we also have that
$$\E^{\Q^*}\big[\int_{0}^{T}|q_s^*|^2\dif s\big]=2\E[M^*_T \ln M^*_T]<+\infty,$$
\begin{align}\label{4a}
\E^{\Q^*}\big[\sup_{t\in\T} A_\beta(t)\big]\leq C_p+\frac{1}{2p} \E^{\Q^*}\Big[\int_0^T|q_s^*|^2 \dif s\Big]<+\infty,
\end{align}
and
\begin{align}\label{4aa}
\E^{\Q^*}\Big[\sup_{t\in\T} \overline{A}_\beta(t)\Big]\leq C_p+\frac{1}{2p} \E^{\Q^*}\Big[\int_0^T|q_s^*|^2 \dif s\Big]<+\infty.
\end{align}
Thus, it suffices to show that
\begin{align}\label{320}
\E^{\Q^*}\Big[\int_0^T|f(s,0,q_s^*)|\dif s\Big]<+\infty.
\end{align}
To do this, we first show that
\begin{align}\label{321}
\E^{\Q^*}\Big[\int_0^T|f(s,Y_s,q_s^*)|\dif s\Big]<+\infty.
\end{align}
In fact, in light of $f(t,y,q)\geq-\overline{\alpha}_t-\beta|y|+\frac{1}{2\gamma}|q|^2$, we deduce
$$
f^-(s,Y_s,q_s^*)\leq\overline{\alpha}_s+\beta|Y_s|,
$$
so, in light of (\ref{4a}) and (\ref{4aa}),
\begin{align}\label{eq:3.27}
\E^{\Q^*}\Big[\int_0^T f^-(s,Y_s,q_s^*)\dif s\Big]\leq\E^{\Q^*}\Big[\int_0^T (\overline{\alpha}_s+\beta|Y_s|)\dif s\Big]<+\infty.
\end{align}
For each $n\geq 1$, we define the following stopping time:
$$
\tau_n:=\inf \Big\{t \in[0, T]: \int_0^t|Z_s|^2 \dif s\geq n\Big\} \wedge T.
$$
By BSDE (\ref{BSDE1.1}) and \eqref{fdy3.1} we have for each $n\geq1$,
\begin{align*}
Y_0=Y_{\tau_n}+\int_0^{\tau_n} f(s,Y_s,q_s^*)\dif s+\int_0^{\tau_n} Z_s\cdot\dif B_s^{q^*}=\E^{\Q^*}\Big[Y_{\tau_n}+\int_0^{\tau_n} f(s,Y_s,q_s^*)\dif s\Big].
\end{align*}
It then follows that
\begin{align*}
Y_0\geq\E^{\Q^*}\Big[-|Y_{\tau_n}|+\int_0^{\tau_n} f^+(s,Y_s,q_s^*)\dif s-\int_0^{\tau_n} f^-(s,Y_s,q_s^*)\dif s\Big].
\end{align*}
Thus, in light of (\ref{4a}), (\ref{4aa}) and \eqref{eq:3.27}, we have for each $n\geq 1$,
\begin{align}\label{eq:3.28}
\E^{\Q^*}\Big[\int_0^{\tau_n} f^+(s,Y_s,q_s^*)\dif s\Big]\leq Y_0+\E^{\Q^*}\Big[\sup_{t\in\T}|Y_t|\Big]+\E^{\Q^*}\Big[\int_0^T f^-(s,Y_s,q_s^*)\dif s\Big]<+\infty.
\end{align}
Combining \eqref{eq:3.27} and \eqref{eq:3.28} along with Fatou's lemma yields (\ref{321}). Furthermore, in light of $f(t,y,q)\geq-\overline{\alpha}_t-\beta|y|+\frac{1}{2\gamma}|q|^2$, we have
$$f(s,0,q_s^*)\geq-\overline{\alpha}_s, ~~s\in[0,T].$$
It then follows from the last inequality and (\ref{x}) that
\begin{align*}
-\overline{\alpha}_s\leq f(s,0,q_s^*)\leq f(s,Y_s,q_s^*)+\beta|Y_s|+\overline{\alpha}_s, ~~s\in\T,
\end{align*}
which implies
$$
|f(s,0,q_s^*)|\leq|f(s,Y_s,q_s^*)|+\beta|Y_s|+\overline{\alpha}_s.
$$
Then, (\ref{320}) follows immediately from the last inequality, (\ref{321}), (\ref{4a}) and (\ref{4aa}). Thus, we have proved the desired assertion that $q^*\in\A$ is optimal for the stochastic control problem, i.e., $Y_\cdot=Y^{q^*}_\cdot=\essinf_{q\in\A}Y^{q}_\cdot$, which naturally yields the uniqueness result. The proof is complete.\hfill\framebox

\section{The uniqueness for the case of $g$ being convex in $(y, z)$}
\label{section 4}
\setcounter{equation}{0}

This section is devoted to the uniqueness for the unbounded solution of BSDE {\rm{($\ref{BSDE1.1}$)}} with generator $g$ being convex in the state variables $(y, z)$.

\subsection{Statement of the main result}

We use the following assumption on $g$.\vspace{0.2cm}

\noindent\textbf{(UN3)} $\as$, for each $(y,z)\in\R\times\R^d, ~g(t,\cdot,\cdot)$ is convex.\vspace{0.2cm}

Suppose that the generator $g$ satisfies (UN3), the Legendre-Fenchel transformation of $g$ can be defined as follows:
\begin{align}\label{fdyc}
f(\omega,t,r,q):=\sup_{(y,z)\in\R\times\R^d}\big(ry+qz-g(\omega,t,y,z)\big),~~~\forall (\omega,t,r,q)\in\Omega\times[0,T]\times\R\times\R^{1\times d}.
\end{align}
The $f$ is a function valued in $\R\cup\{+\infty\}$.
We define the admissible control set:
\begin{align*}
\A:=\bigg\{&(r_t,q_t)_{t\in[0,T]}~~\text{is a pair of } (\F_t)\text{-progressively measurable~}  \R\times\R^{1\times d} \text{-valued processes}:\\ &-\beta\leq r_\cdot\leq\beta,~~\ps,~\int_{0}^{T}|q_s|^2\dif s<+\infty,~~\E^{\Q}\Big[\int_{0}^{T}|q_s|^2\dif s\Big]<+\infty,\\
&\E^{\Q}\Big[|\xi|+\int_{0}^{T}|f(s,r_s,q_s)|\dif s\Big]<+\infty, ~\big(M_t\big)_{t\in[0,T]}~~\text{is a uniformly integrable martingale}\\
&\text{with} ~~M_t:=\exp\Big(\int_{0}^{t}q_s\dif B_s-\frac{1}{2}\int_{0}^{t}|q_s|^2\dif s\Big), ~t\in\T, ~\text{and}~\frac{\dif \Q}{\dif \mathbb{P}}:=M_{T}\bigg\}.
\end{align*}

The following Theorem \ref{unthc} gives a uniqueness result for the unbounded solution of BSDE (\ref{BSDE1.1}) under assumption (UN3).

\begin{theorem}\label{unthc}
Suppose that $\xi$ is a terminal value, the generator $g$ is continuous in $(y,z)$ and satisfies {\rm{(EX)}}, $\E[\exp(pe^{\beta T}(\xi^-+\int_0^T\overline{\alpha}_t\dif t))]<+\infty$ for some $p>\gamma$, and $\E[\exp(\overline{p}e^{\beta T}( \xi^++\int_0^T\underline{\alpha}_t\dif t))]<+\infty$ for any $\overline{p}>\gamma$. If $g$ also satisfies {\rm{(UN2)}} and {\rm{(UN3)}}, then BSDE {\rm{($\ref{BSDE1.1}$)}} admits a unique adapted solution $(Y_\cdot,Z_\cdot)$ satisfying {\rm{(\ref{p1})}} and {\rm{(\ref{p2})}} for any $\overline{p}>\gamma$. Moreover, we have an explicit expression of $Y_\cdot$:
\begin{align*}
Y_t=\essinf_{(r,q)\in\A}Y_t^{r,q}, ~~t\in[0,T],\vspace{-0.2cm}
\end{align*}
where
\begin{align*}
 Y_t^{r,q}:=\E^\Q\Big[e^{-\int_t^Tr_s\dif s} \xi+\int_t^Te^{-\int_t^sr_u\dif u} f(s,r_s,q_s)\dif s\big| \F_t\Big]
\end{align*}
and the function $f$ is defined in \eqref{fdyc}.\vspace{0.2cm}
\end{theorem}

\begin{remark}\label{rmk-lip}
{\rm{With respect to Theorem \ref{unthc}, we make the following two comments.
\begin{itemize}
\item[(i)] If the generator $g$ satisfies (1) of (EX) and (UN2), then $\as$, we have
\begin{align*}
|g(t,y,z)|\leq\underline{\alpha}_t+\overline{\alpha}_t+\beta|y|+\frac{\gamma}{2}|z|^2, ~(y,z)\in\R_+\times\R^d,
\end{align*}
which implies that $g$ satisfies the linear growth condition in $y$ when $y\geq0$.
\item[(ii)] If the generator $g$ satisfies (2) of (EX) and (UN2), then $\as$, we have
\begin{align*}
g(t,y,z)\leq\overline{\alpha}_t+\beta|y|+\frac{\gamma}{2}|z|^2, ~(y,z)\in\R\times\R^d,
\end{align*}
which along with (UN3) gives that $\as$,
\begin{align}\label{lip}
|g(t,y_1,z)-g(t,y_2,z)|\leq \beta|y_1-y_2|, ~~~(y_1, y_2,z)\in\R\times\R\times\R^d.
\end{align}
In fact, $\as$, for each $(y_1, y_2,z)\in\R\times\R\times\R^d$ and $\theta\in(0,1)$, we have
\begin{align*}
g(t,y_1,z)&=g\Big(t,\theta y_2+(1-\theta)\frac{y_1-\theta y_2}{1-\theta},\theta z+(1-\theta)z\Big)\\
&\leq\theta g(t,y_2,z)+(1-\theta)g\Big(t,\frac{y_1-\theta y_2}{1-\theta},z\Big)\\
&\leq\theta g(t,y_2,z)+(1-\theta)\Big(\overline{\alpha}_t+\beta\Big|\frac{y_1-\theta y_2}{1-\theta}\Big|+\frac{\gamma}{2}|z|^2\Big).
\end{align*}
Letting $\theta\rightarrow1$ in the last inequality yields that $g(t,y_1,z)-g(t,y_2,z)\leq \beta|y_1-y_2|$. Then, the desired inequality (\ref{lip}) follows by exchanging the positions of $y_1$ and $y_2$. Consequently, the assumptions in Theorem \ref{unthc} are strictly stronger than the uniformly Lipschitz continuity condition in the state variable $y$ used in \cite{Delbaen 2011}. However, we can get a more accurate conclusion in this situation, that is, $Y_\cdot$ has a more explicit expression.
\end{itemize}}}
\end{remark}

\subsection{Proof of the main result}

In this subsection, we will give the proof of Theorem \ref{unthc}.\vspace{0.2cm}

\noindent\textbf{Proof of Theorem \ref{unthc}.}
The existence result has been given in Proposition \ref{cun1}. Now, we show the uniqueness part. According to assumptions {\rm{(EX)}} and {\rm{(UN3)}}, we can deduce that the function $f$ defined in (\ref{fdyc}) verifies the following properties: $\as$,
\begin{itemize}
\item[(1)] $\forall (r,q)\in\R\times{\R^{1\times d}}, ~f(t,r,q)\geq-\overline{\alpha}_t+\frac{1}{2\gamma}|q|^2$;
\item[(2)] If $f(t,r,q)<+\infty$ for some $(r,q)\in\R\times\R^{1\times d}$, then $r\in[-\beta,\beta]$;
\item[(3)] $f(t,\cdot,\cdot)$ is convex.
\end{itemize}
Since (3) is trivial, we need only to show (1) and (2) here. In fact, in light of (3) of (EX), we have $\as$, for each $z\in\R^d, ~|g(t,0,z)|\leq\overline{\alpha}_t+\frac{\gamma}{2}|z|^2$. It then follows that $\as$, for each $(r,q)\in\R\times\R^{1\times d}$,
\begin{align*}
f(t,r,q)&=\sup_{(y,z)\in\R\times\R^d}\big(ry+qz-g(t,y,z)\big)\geq\sup_{z\in\R^d}
\big(qz-g(t,0,z)\big)\\
&\geq\sup_{z\in\R^d}\Big(qz-\overline{\alpha}_t-\frac{\gamma}{2}|z|^2\Big)
=-\overline{\alpha}_t+\frac{1}{2\gamma}|q|^2.
\end{align*}
Then, (1) is true. Furthermore, thanks to (2) of (EX) and {\rm{(UN2)}}, we have $\as$, for each $(r,q)\in\R\times\R^{1\times d}$,
\begin{align*}
f(t,r,q)&=\sup_{(y,z)\in\R\times\R^d}\big(ry+qz-g(t,y,z)\big)\\
&\geq\max\Big\{\sup_{y\leq0}\big(ry-g(t,y,0)\big), ~\sup_{y>0}\big(ry-g(t,y,0)\big)\Big\}\\
&\geq\max\Big\{\sup_{y\leq0}\big(ry-\overline{\alpha}_t+\beta y\big), ~\sup_{y>0}\big(ry-\overline{\alpha}_t-\beta y\big)\Big\}\\
&=-\overline{\alpha}_t+\max\Big\{\sup_{y\leq0}\big(r+\beta \big)y, ~\sup_{y>0}\big(r-\beta\big)y\Big\}.\vspace{0.1cm}
\end{align*}
Hence, if $f(t,r,q)<+\infty$, then $r+\beta\geq0$ and $r-\beta\leq0$, i.e., $-\beta\leq r\leq\beta$. Thus, (2) is true.\vspace{0.2cm}

Now, let $(r, q)$ be in $\A$, if this set is not empty. We set $B_t^q:=B_t-\int_0^t q_s\dif s,~t\in\T$. According to Girsanov's theorem, we can conclude that $(B_t^q)_{t\in\T}$ is a standard $d$-dimensional Brownian motion under the probability measure $\Q$. Let us consider the following BSDE:
\begin{align}\label{bsde c}
Y_t^{r,q}=\xi+\int_{t}^{T}\big(f(s,r_s,q_s)-r_sY_s^{r,q}\big)\dif s+\int_{t}^{T}Z_s^{r,q}\cdot\dif B_s^q,~~~~t\in\T.
\end{align}
In light of the definition of $\A$, by Theorem 6.2 in \cite{Briand 2003} we know that BSDE {\rm{(\ref{bsde c})}} admits a unique adapted solution $(Y_t^{r,q},Z_t^{r,q})_{t\in\T}$ such that $Y_\cdot^{r,q}$ belongs to the class $(D)$ under $\Q$. Hence, the stochastic control problem consists in minimizing $Y_\cdot^{r,q}$ among all the admissible controls $(r, q)\in\A$. In what follows, let $(Y_t,Z_t)_{t\in\T}$ be an adapted solution of BSDE {\rm{(\ref{BSDE1.1})}} satisfying (\ref{p1}) and (\ref{p2}) for any $\overline{p}>\gamma$. We divide our proof into the following three steps.\vspace{0.2cm}

\textbf{First step.} Let us start with showing $Y_\cdot\leq Y_\cdot^{r,q}$ for any $(r,q)\in\A$. Firstly, we have $Y_T=Y^{r,q}_T=\xi$. For each $t\in[0,T]$ and each integer $m\geq1$, denote the following stopping time
$$\tau_m^t:=\inf\Big\{s\geq t, \int_t^s|Z_u|^2\dif u+\int_t^s|Z^{r,q}_u|^2\dif u+\int_t^s|q_u|^2\dif u\geq m\Big\}\wedge T$$
with the convention $\inf\emptyset=+\infty$. By applying $\rm{It\hat{o}}$'s formula to $e^{-\int_0^t r_s\dif s} Y_t^{r,q}$ we obtain
\begin{align*}
e^{-\int_0^t r_s\dif s} Y_t^{r,q}=e^{-\int_0^{\tau_m^t} r_s\dif s} Y_{\tau_m^t}^{r,q}+\int_t^{\tau_m^t}e^{-\int_0^s r_u\dif u} f(s,r_s,q_s)\dif s+\int_t^{\tau_m^t}e^{-\int_0^s r_u\dif u} Z_s^{r,q}\cdot\dif B_s^q.
\end{align*}
It then follows that
\begin{align*}
Y_t^{r,q}=\E^{\Q}\Big[e^{-\int_t^{\tau_m^t} r_s\dif s} Y_{\tau_m^t}^{r,q}+\int_t^{\tau_m^t}e^{-\int_t^s r_u\dif u} f(s,r_s,q_s)\dif s\big|\F_t\Big].
\end{align*}
Since $e^{-\int_t^{\tau_m^t} r_s\dif s}\leq e^{\int_t^{\tau_m^t} \beta\dif s}\leq e^{\beta T}$ and $Y_\cdot^{r,q}$ belongs to the class $(D)$ under $\Q$, by letting $m\rightarrow\infty$ in the last inequality, we have
\begin{align}\label{yrq}
Y_t^{r,q}=\E^{\Q}\Big[e^{-\int_t^T r_s\dif s} \xi+\int_t^Te^{-\int_t^s r_u\dif u} f(s,r_s,q_s)\dif s\big|\F_t\Big].
\end{align}
Similarly, $\rm{It\hat{o}}$'s formula applied to $e^{-\int_0^t r_s\dif s} Y_t$ gives, in light of (\ref{fdyc}),
\begin{align*}
e^{-\int_0^t r_s\dif s} Y_t&=e^{-\int_0^{\tau_m^t} r_s\dif s} Y_{\tau_m^t}+\int_t^{\tau_m^t}e^{-\int_0^s r_u\dif u} (-g(s,Y_s,Z_s)+r_sY_s+q_sZ_s)\dif s\\
&~~~~~+\int_t^{\tau_m^t}e^{-\int_0^s r_u\dif u} Z_s\cdot\dif B_s^q\\
&\leq e^{-\int_0^{\tau_m^t} r_s\dif s} Y_{\tau_m^t}+\int_t^{\tau_m^t}e^{-\int_0^s r_u\dif u} f(s,r_s,q_s)\dif s+\int_t^{\tau_m^t}e^{-\int_0^s r_u\dif u} Z_s\cdot\dif B_s^q.
\end{align*}
So, we have
\begin{align}\label{e}
Y_t\leq\E^{\Q}\Big[e^{-\int_t^{\tau_m^t} r_s\dif s} Y_{\tau_m^t}+\int_t^{\tau_m^t}e^{-\int_t^s r_u\dif u} f(s,r_s,q_s)\dif s\big|\F_t\Big].
\end{align}
In light of (\ref{p1}) and (\ref{p2}), Fenchel's inequality gives, picking $\overline{p}=p$,
\begin{align*}
\E^{\Q}\Big[\sup_{s\in\T} Y_s^-+\sup_{s\in\T} Y_s^+\Big]&=\E\Big[M_T \sup_{s\in\T} Y_s^-\Big]+\E\Big[M_T \sup_{s\in\T} Y_s^+\Big]\\
&\leq \E\Big[\exp\Big(p\sup_{s\in\T} Y_s^+\Big)\Big]+\E\Big[\exp\Big(p\sup_{s\in\T} Y_s^-\Big)\Big]\\
&~~~~+\frac{2}{p}\E[M_T(\ln M_T-\ln p-1)]\\
&\leq C_p+\frac{1}{p}\E^{\Q}\Big[\int_{0}^{T}|q_s|^2\dif s\Big]<+\infty.
\end{align*}
Then, noticing that $0\leq e^{-\int_t^s r_u\dif u} \leq e^{\beta T}$ for each $s\in [t,\tau_m^t]$, by letting $m\rightarrow\infty$ and applying the Lebesgue's dominated convergence theorem in (\ref{e}), we obtain
\begin{align}\label{y}
Y_t\leq\E^{\Q}\Big[e^{-\int_t^T r_s\dif s} \xi+\int_t^Te^{-\int_t^s r_u\dif u} f(s,r_s,q_s)\dif s\big|\F_t\Big].
\end{align}
Combining (\ref{yrq}) and (\ref{y}), we deduce that for each $t\in\T, ~\ps,~Y_t\leq Y_t^{r,q}$.\vspace{0.2cm}

\noindent\textbf{Second step.} Set $(r_s^*, q_s^*)\in\partial_{(y,z)}g(s,Y_s,Z_s)$. Then we have
$$
f(s,r_s^*,q_s^*)=r_s^*Y_s+q_s^*Z_s-g(s,Y_s,Z_s), ~~s\in[0,T].
$$
Consequently, if $(r^*, q^*)\in\A$, then by BSDEs (\ref{BSDE1.1}) and (\ref{yrq}) along with the above argument we can conclude that $Y_\cdot=Y_\cdot^{r^*, q^*}$.\vspace{0.2cm}

\noindent\textbf{Third step.} We conclude the proof by verifying $(r^*, q^*)\in\A$.
Since $(r_s^*, q_s^*)\in\partial_{(y,z)}g(s,Y_s,Z_s)$, it is clear that $f(s,r_s^*, q_s^*)<+\infty$, then we can deduce from (2) that $-\beta\leq r_s^*\leq\beta$ and, in light of (1),
\begin{align*}
g(s,Y_s,Z_s)&=r_s^*Y_s+q_s^*Z_s-f(s,r_s^*,q_s^*)\\
&\leq r_s^*Y_s+\frac{1}{2}\Big(\frac{1}{2\gamma}|q_s^*|^2+2\gamma|Z_s|^2\Big)
+\overline{\alpha}_s-\frac{1}{2\gamma}|q_s^*|^2\\
&= r_s^*Y_s+\gamma|Z_s|^2+\overline{\alpha}_s-\frac{1}{4\gamma}|q_s^*|^2;\\
\frac{1}{4\gamma}|q_s^*|^2&\leq-g(s,Y_s,Z_s)+\beta |Y_s|+\gamma|Z_s|^2+\overline{\alpha}_s.
\end{align*}
It follows from the last inequality that $\ps,~\int_0^T|q_s^*|^2\dif s<+\infty$. For each $n\geq1$, we define
$$
M^*_t:=\exp\Big(\int_{0}^{t}q_s^*\dif B_s-\frac{1}{2}\int_{0}^{t}|q_s^*|^2\dif s\Big), ~~t\in\T,
$$
$$
\tau_n:=\inf \Big\{t \in[0, T]: \int_0^t|q_s^*|^2 \dif s+\int_0^t|Z_s|^2 \dif s\geq n\Big\} \wedge T, \quad \frac{\dif \Q_n^*}{\dif \mathbb{P}}:=M^*_{\tau_n}.\vspace{0.1cm}
$$
Set $B_t^{q^*}:=B_t-\int_0^t q^*_s\dif s,~t\in\T$. Thanks to Girsanov's theorem, $(B_t^{q^*})_{t\in [0,\tau_n]}$ is a standard $d$-dimensional Brownian motion under $\Q^*_n$. Let us further show that $\big(M^*_t\big)_{t\in[0,T]}$ is a uniformly integrable martingale. To do this, we need to prove the following Lemma.

\begin{lemma}\label{lem c}
$(M^*_{\tau_n})_n$ is uniformly integrable.
\end{lemma}

\noindent\textbf{Proof of Lemma \ref{lem c}.}
According to (\ref{p1}), (\ref{p2}) and Fenchel's inequality, we have
\begin{align}\label{a}
\E^{\Q_n^{*}}\Big[\sup_{t\in\T} A_\beta(t)\Big]=\E\Big[M^*_{\tau_n} \sup_{t\in\T} A_\beta(t)\Big] \leq C_p+\frac{1}{2p} \E^{\Q_n^*}\Big[\int_0^{\tau_n}|q_s^*|^2 \dif s\Big],
\end{align}
and for each $\overline{p}>\gamma$, in the same manner,
\begin{align}\label{aa}
\E^{\Q_n^{*}}\Big[\sup_{t\in\T} \overline{A}_\beta(t)\Big]\leq C_{\overline{p}}+\frac{1}{2\overline{p}} \E^{\Q_n^*}\Big[\int_0^{\tau_n}|q_s^*|^2 \dif s\Big].
\end{align}
Since $g(s,Y_s,Z_s)=r_s^*Y_s+q_s^*Z_s-f(s,r_s^*,q_s^*)$, and $(M^*_{t\wedge\tau_n})_{t\in[0,T]}$ is a martingale under $\mathbb{P}$, we have
\begin{align}\label{eq:4.10}
Y_{t\wedge\tau_n}=Y_{\tau_n}+\int_{t\wedge\tau_n}^{\tau_n}\big(-r_s^*Y_s
+f(s,r_s^*,q_s^*)\big)\dif s+\int_{t\wedge\tau_n}^{\tau_n} Z_s \cdot \dif B_s^{q^*}, ~~t\in\T.
\end{align}
Applying $\rm{It\hat{o}}$'s formula to $e^{\beta s} Y_s$, we obtain
\begin{align*}
\dif \big(e^{\beta s} Y_s\big)=e^{\beta s}\big((\beta+r_s^*) Y_s-f(s,r_s^*,q_s^*)\big)\dif s-e^{\beta s}Z_s\cdot\dif B_s^{q^*}, ~~s\in[0,\tau_n].
\end{align*}
Moreover, in light of $f(t,r,q)\geq-\overline{\alpha}_t+\frac{1}{2\gamma}|q|^2$, $Y_{\tau_n}\geq-Y_{\tau_n}^-$ and $0\leq \beta+r_s^*\leq 2\beta$, we deduce that
\begin{align*}
Y_0&=e^{\beta \tau_n} Y_{\tau_n}+\int_0^{\tau_n}e^{\beta s}\big(f(s,r_s^*,q_s^*)-(\beta+r_s^*) Y_s\big)\dif s+\int_0^{\tau_n}e^{\beta s}Z_s\cdot\dif B_s^{q^*}\\
&\geq-e^{\beta \tau_n} Y_{\tau_n}^-+\int_0^{\tau_n}e^{\beta s}\big(-\overline{\alpha}_s-2\beta Y_s^++\frac{1}{2\gamma}|q_s^*|^2\big)\dif s+\int_0^{\tau_n}e^{\beta s}Z_s\cdot\dif B_s^{q^*}.
\end{align*}
It then follows from (\ref{a}) and (\ref{aa}) that for each $\overline{p}>\gamma$, similar to (\ref{tt}),
\begin{align*}
\notag Y_0&\geq \E^{\Q_n^*}\Big[-e^{\beta \tau_n} Y_{\tau_n}^-+\int_0^{\tau_n}e^{\beta s}\big(-\overline{\alpha}_s-2\beta Y_s^++\frac{1}{2\gamma}|q_s^*|^2\big)\dif s\Big]\\
&\geq-C_p-2\beta TC_{\overline{p}}+\frac{1}{2}\Big(\frac{1}{\gamma}-\frac{1}{p}-\frac{2\beta T}{\overline{p}}\Big)\E^{\Q_n^*}\Big[\int_0^{\tau_n}|q_s^*|^2\dif s\Big].
\end{align*}
Set $\overline{p}>\frac{2\beta Tp\gamma }{p-\gamma}$. Then $\frac{1}{2}\Big(\frac{1}{\gamma}-\frac{1}{p}-\frac{2\beta T}{\overline{p}}\Big)>0.$ Hence,
\begin{align}\label{g}
2\E[M^*_{\tau_n}\ln M^*_{\tau_n}]=2\E^{\Q_n^{*}}[\ln M^*_{\tau_n}]=\E^{\Q_n^{*}}\Big[\int_0^{\tau_n}|q_s^*|^2\dif s\Big]<C_{p,\overline{p},T,\beta,\gamma},
\end{align}
where $C_{p,\overline{p},T,\beta,\gamma}$ is a positive constant independent of $n$. Hence, we obtain $\sup_n\E[M^*_{\tau_n}\ln M^*_{\tau_n}]<+\infty$. Then we conclude the proof of the lemma by using the de La Vall$\acute{\rm{e}}$e Poussin lemma.
\hfill\framebox

\vspace{0.2cm}

Applying Fatou's lemma to (\ref{g}), we obtain
\begin{align}\label{eq:4.12}
2\E[M^*_T\ln M^*_T]=\E^{\Q^*}\Big[\int_0^T|q_s^*|^2\dif s\Big]\leq\liminf_n\E^{\Q_n^* }\Big[\int_0^{\tau_n}|q_s^*|^2\dif s\Big]<+\infty.
\end{align}
Thus, $\big(M^*_t\big)_{t\in[0,T]}$ is a uniformly integrable martingale and defines a probability measure $\Q^*\thicksim \mathbb{P}$ via
$$\frac{\dif \Q^*}{\dif \mathbb{P}}:=M^*_T.$$
According to Girsanov's theorem, $(B_t^{q^*})_{t\in\T}$ is a standard $d$-dimensional Brownian motion under $\Q^*$. Finally, it remains to show that $\E^{\Q^*}\big[\int_0^T|f(s,r_s^*,q_s^*)|\dif s\big]<+\infty$. In fact, in light of \eqref{eq:4.12} and the fact that $f(t,r,q)\geq-\overline{\alpha}_t+\frac{1}{2\gamma}|q|^2$, by Fenchel's inequality we deduce that
$$
f^-(s,r_s^*,q_s^*)\leq\overline{\alpha}_s, ~~s\in\T,
$$
and then
\begin{align}\label{eq:4.13}
\E^{\Q^*}\Big[\int_0^T f^-(s,r_s^*,q_s^*)\dif s\Big]\leq\E^{\Q^*}\Big[\int_0^T \overline{\alpha}_s\dif s\Big]<+\infty.
\end{align}
In light of (\ref{p1}), (\ref{p2}) and \eqref{eq:4.12}, Fenchel's inequality gives, picking $\overline{p}=p$,
\begin{align}\label{eq:4.14}
\E^{\Q^*}\Big[\sup_{s\in\T} Y_s^-+\sup_{s\in\T} Y_s^+\Big]&=\E\Big[M_T^* \sup_{s\in\T} Y_s^-\Big]+\E\Big[M_T^* \sup_{s\in\T} Y_s^+\Big]\nonumber\\
&\leq \E\Big[\exp\Big(p\sup_{s\in\T} Y_s^+\Big)\Big]+\E\Big[\exp\Big(p\sup_{s\in\T} Y_s^-\Big)\Big]\nonumber\\
&~~~~+\frac{2}{p}\E[M_T^*(\ln M_T^*-\ln p-1)]\\
&<+\infty.\nonumber\vspace{-0.2cm}
\end{align}
Furthermore, it follows from \eqref{eq:4.10} and the fact $-\beta\leq r_\cdot^*\leq\beta$ that for each $n\geq 1$,
\begin{align*}
Y_0&=\E^{\Q^*}\Big[e^{-\int_0^{\tau_n} r_s^*\dif s}Y_{\tau_n}+\int_0^{\tau_n} e^{-\int_0^s r_u^*\dif u} f(s,r_s^*,q_s^*)\dif s\Big]\\
&\geq\E^{\Q^*}\Big[-e^{\beta T} |Y_{\tau_n}|+ e^{-\beta T} \int_0^{\tau_n} f^+(s,r_s^*,q_s^*)\dif s- e^{\beta T} \int_0^{\tau_n} f^-(s,r_s^*,q_s^*)\dif s\Big].
\end{align*}
By the last inequality along with \eqref{eq:4.13} and \eqref{eq:4.14} we obtain that for each $n\geq 1$,
\begin{align}\label{eq:4.15}
&\E^{\Q^*}\Big[\int_0^{\tau_n} f^+(s,r_s^*,q_s^*)\dif s\Big]\nonumber \\
&\ \ \ \leq e^{\beta T}Y_0+e^{2\beta T}\E^{\Q^*}\Big[\sup_{t\in\T}|Y_t|\Big]+e^{2\beta T}\E^{\Q^*}\Big[\int_0^T f^-(s,r_s^*,q_s^*)\dif s\Big]<+\infty.
\end{align}
Combining \eqref{eq:4.13} and \eqref{eq:4.15} along with Fatou's lemma yields that $\E^{\Q^*}\big[\int_0^T|f(s,r_s^*,q_s^*)|\dif s\big]<+\infty$. Thus, we show that $(r^*, q^*)\in\A$ is optimal, i.e., for each $t\in\T, ~\ps,~Y_t=Y^{r^*, q^*}_t=\essinf_{(r, q)\in\A}Y^{r,q}_t$, which naturally yields the uniqueness result. The proof of Theorem \ref{unthc} is complete.\hfill\framebox

\section{Comparison theorem}
\label{section 5}
\setcounter{equation}{0}

In this section, we aim at establishing a general comparison theorem for unbounded solutions of BSDE {\rm{($\ref{BSDE1.1}$)}} under the assumptions used in Section \ref{section 3}. For this purpose, we first apply Proposition 5 in \cite{Briand-Hu 2006} to obtain the following result.

\begin{proposition}\label{cth1}
Let $\xi$ and $\xi'$ be two terminal values, $g$ and $g'$ be two generators, and $(Y_\cdot,Z_\cdot)$ and $(Y'_\cdot,Z'_\cdot)$ be, respectively, an adapted solution of BSDE $(\xi, T, g)$ and BSDE $(\xi', T, g')$. We assume that $\ps,~\xi\leq\xi'$ and that $g'$ satisfies, $\as$, for some constants $\mu$ and $\lambda$,
\begin{align*}
&-\sgn(y-y')(g'(t,y,z)-g'(t,y',z))\leq\mu|y-y'|;\\
&|g'(t,y,z)-g'(t,y,z')|\leq\lambda|z-z'|.
\end{align*}
If $(Y_\cdot-Y'_\cdot)^+\in\s$ and $\as$,
$$
-g(t,Y_t,Z_t)\leq -g'(t,Y_t,Z_t),
$$
then for each $t\in[0,T]$, $\ps$, we have $Y_t\leq Y'_t$.
\end{proposition}

The following comparison theorem for unbounded solutions of BSDE {\rm{($\ref{BSDE1.1}$)}} is the main result of this section.

\begin{theorem}\label{cth}
Let $\xi$ and $\xi'$ be two terminal values, $g$ and $g'$ be two generators such that $g$ is continuous in $(y,z)$ and satisfies assumptions {\rm{(EX)}}, {\rm{(UN1)}} and {\rm{(UN2)}}, and $(Y_\cdot,Z_\cdot)$ and $(Y'_\cdot,Z'_\cdot)$  be, respectively, an adapted solution of BSDE $(\xi, T, g)$ and BSDE $(\xi', T, g')$ such that both {\rm{(\ref{p1})}} and {\rm{(\ref{p2})}} hold for some $p>\gamma$ and any $\overline{p}>\gamma$ and $Y'_\cdot\in\s$. We further assume that $\as,$ for each $(y,z)\in\R\times\R^d$,
\begin{align}\label{bj}
g(\omega,t,y,z)\geq-\overline{\alpha}_t(\omega)-\varphi(|y|)-\gamma|z|.
\end{align}
If $\ps,~\xi'\leq\xi$ and $\as$,
$$
-g'(t,Y'_t,Z'_t)\leq -g(t,Y'_t,Z'_t),
$$
then for each $t\in[0,T]$, $\ps$, we have $Y'_t\leq Y_t$.
\end{theorem}

\begin{remark}\label{rmk-cth}
{\rm{With respect to Theorem \ref{cth}, we make the following three comments.
\begin{itemize}
\item[(i)] It should be noted that the inequality (\ref{bj}) proposed in Theorem \ref{cth} does not impose too many restrictions. For example, we consider the case that $g$ only depends on $z$. If $g(z)$ is convex, then there must exist a constant $c$ such that
$$
g(z)\geq-c|z|-c,
$$
which naturally yields (\ref{bj}).
\item[(ii)] Thanks to (ii) of Theorem \ref{unth3.1}, the conclusions of Theorem \ref{cth} still hold without assuming (UN2) if the (\ref{p2}) appearing in Theorem \ref{cth} is replaced by $Y_\cdot\leq K$ and $Y'_\cdot \leq K$ for some constant $K>0$. The proof is similar to the proof of Theorem \ref{cth}.
\item[(iii)] Thanks to Theorem \ref{unthc}, the conclusions of Theorem \ref{cth} still hold if the assumption (UN1) in Theorem \ref{cth} is replaced by (UN3). In fact, for each integer $n\geq1$, we can define $g_n$ as follows:
\begin{align*}
g_n(\omega,t,y,z):=\inf_{(u,v)\in\R\times\R^d}\big\{g(\omega,t,u,v)+(n+\beta)|y-u|
+(n+\gamma)|z-v|\big\}.
\end{align*}
More specifically, in light of (ii) of Remark \ref{rmk-lip}, we know that $g$ is Lipschitz continuous in $y$ under assumptions (EX) and (UN2), which implies that $\as,$ for each $(y,z)\in\R\times\R^d$,
\begin{align}\label{q1}
g(t,0,z)-g(t,y,z)\leq \beta|y|.
\end{align}
It follows from the inequality (\ref{bj}) that $g(t,0,z)\geq-\overline{\alpha}_t-\gamma|z|$, which combined with (\ref{q1}) yields
\begin{align*}
g(t,y,z)\geq-\overline{\alpha}_t-\beta|y|-\gamma|z|.
\end{align*}
Hence, we can deduce that $g_n$ is well defined, $g_n\uparrow g$, and $g_n$ is $(n+\beta)$-Lipschitz continuous in $y$ and $(n+\gamma)$-Lipschitz continuous in $z$. Thus, the subsequent proof of Theorem \ref{cth} can proceed smoothly.
\end{itemize}}}
\end{remark}

In the sequel, we present the proof of Theorem \ref{cth}.\vspace{0.2cm}

\noindent\textbf{Proof of Theorem \ref{cth}.}
For each integer $n\geq1$, let
$$
g_n(\omega,t,y,z):=\inf_{v\in\R^d}\big\{g(\omega,t,y,v)+(n+\gamma)|z-v|\big\}, ~~~(\omega,t,y,z)\in\Omega\times\T\times\R\times\R^d.
$$
Thanks to $g(t,y,\cdot)$ being convex and (\ref{bj}), $g_n$ is well defined. In light of $g$ satisfying {\rm{(EX)}} and {\rm{(UN1)}}, it can be easily checked that $g_n$ satisfies {\rm{(EX)}} and (1) of {\rm{(UN1)}}. By the proof of Theorem 1 in \cite{Fan-Jiang-Tian 2011}, it is not hard to verify that $\{g_n\}_{n\geq1}$ is nondecreasing and converges to $g$, and $g_n$ satisfies that $\as$, for each $(y,z_1,z_2)\in\R\times\R^d\times\R^d$,
\begin{align*}
|g_n(t,y,z_1)-g_n(t,y,z_2)|\leq(n+\gamma)|z_1-z_2|.
\end{align*}

Moreover, in light of $g_n\leq g$, we have
$$-g'(t,Y'_t,Z'_t)\leq -g(t,Y'_t,Z'_t)\leq-g_n(t,Y'_t,Z'_t).$$
Thanks to Theorem 4.2 in \cite{Briand 2003} and Proposition \ref{cun1}, we can conclude that BSDE $(\xi,T,g_n)$ admits a unique adapted solution $(Y^n_t,Z^n_t)_{t\in\T}$ such that (\ref{p1}) and (\ref{p2}) hold with $Y^n_\cdot$ instead of $Y_\cdot$ for some $p>\gamma$ and any $\overline{p}>\gamma$. Hence, in light of $Y'_\cdot\in\s$, it is not hard to check that $(Y'_\cdot-Y^n_\cdot)^+\in\s$. Applying Proposition \ref{cth1} gives that for each $n\geq1$, $Y'_\cdot\leq Y^n_\cdot$.

Similarly, since $-g(t,Y_t,Z_t)\leq-g_n(t,Y_t,Z_t)$, we know that for each $n\geq1$, $Y_\cdot\leq Y^n_\cdot$. By virtue of $(g_n)_{n\geq1}$ being nondecreasing, we know that $(Y^n_t)_{n\geq1}$ is non-increasing. Hence, $Y_\cdot\leq Y^n_\cdot\leq Y^1_\cdot$. Set
$$\widehat{Y}_t:=\lim\limits_{n\rightarrow+\infty}Y^n_t,\ \ t\in \T.$$
Since both $Y_\cdot$ and $Y^1_\cdot$ satisfy (\ref{p1}) and (\ref{p2}) for some $p>\gamma$ and any $\overline{p}>\gamma$, proceeding with a localization argument, we obtain that there exists an $(\F_t)$-progressively measurable process $(\widehat{Z}_t)_{t\in\T}$ such that $(\widehat{Y}_\cdot,\widehat{Z}_\cdot)$ is a solution of BSDE $(\xi,T,g)$ satisfying (\ref{p1}) and (\ref{p2}) for some $p>\gamma$ and any $\overline{p}>\gamma$. Furthermore, thanks to (i) of Theorem \ref{unth3.1} along with assumptions (EX), (UN1) and (UN2) of the generator $g$, we know that the adapted solution of BSDE $(\xi,T,g)$ satisfying the above-mentioned conditions is unique, which implies that $\widehat{Y}_\cdot=Y_\cdot$, i.e., $Y_\cdot=\lim_{n\rightarrow+\infty}Y^n_\cdot$. Combining with the fact that $Y'_\cdot\leq Y^n_\cdot$, we can conclude that for each $t\in\T, ~\ps,~Y'_t\leq Y_t$. The proof is then complete.\hfill\framebox

\section{Conclusions}
\label{section 6}
\setcounter{equation}{0}

This paper focuses on the existence and uniqueness of unbounded solutions for quadratic BSDEs, whose generator $g$ may be neither linearly growing nor locally-Lipschitz continuous in the state variable $y$, and is convex with a quadratic growth in the state variable $z$.

Suppose that the terminal value $\xi^-$ admits a certain exponential moment and $\xi^+$ admits all exponential moments or is bounded, and that the generator $g$ is convex with a quadratic growth in the state variable $z$. Under some mild growth conditions on the generator $g$, we obtain the following results: If $g$ satisfies the monotonicity condition with a general growth in the state variable $y$, which is strictly weaker than the uniformly Lipschitz continuity condition used in \cite{Delbaen 2011}, without imposing  the locally-Lipschitz continuity condition of $g$ in $y$ nor the extended convexity condition \eqref{eq:1.3} proposed in \cite{Fan-Hu-Tang 2020}, then the quadratic BSDE admits a unique adapted solution $(Y_t,Z_t)_{t\in\T}$ such that $Y_\cdot$ satisfies certain integrability. If $g$ is jointly convex in the state variables $(y,z)$, which is strictly stronger than the uniformly Lipschitz continuity condition in $y$ when combined with other necessary assumptions, then we can not only have existence and uniqueness of the unbounded solution, but also get an explicit expression of $Y_\cdot$. The Legendre-Fenchel transform of convex functions, Girsanov's theorem, the de La Vall$\acute{\rm{e}}$e Poussin lemma and Fenchel's inequality play important roles in the proof of these results.

In comparison to \cite{Delbaen 2011} and \cite{Delbaen 2015}, we innovatively use the exponential shifting transform with a factor $e^{\beta t}$ or $e^{-\int_0^t r_s\dif s}$ to make full use of the monotonicity condition of the generator $g$ in $y$ so that $g$ is allowed to have a general growth rather than the usual linear growth in $y$, which is one of the key ideas of this paper. Another point that should be stressed is that the time interval partitioning method in \cite{Delbaen 2011} plays a crucial role in establishing the results on the uniform integrability. However, for simplicity of presentation, we do not use this method any longer, and instead lift the integrability of $\xi^+$ so that the proof can proceed smoothly. Moreover, several examples are provided to illustrate that our results are more general than those in \cite{Delbaen 2011}. Finally, it is the first time, to the best of our knowledge, that the corresponding comparison theorems for unbounded solutions of the preceding BSDEs are established under general assumptions.

It is obvious that $(Y_\cdot,Z_\cdot)$ is a solution of BSDE $(\xi, T, g)$ if and only if $(-Y_\cdot,-Z_\cdot)$ is a solution of BSDE $(-\xi, T, \overline{g})$, where $\overline{g}(t,y,z):=-g(t,-y,-z)$, and that $g$ is convex in $z$ (resp. $(y,z)$) if and only if $\overline{g}$ is concave in $z$ (resp. $(y,z)$). Consequently, from this point of view, Theorems \ref{unth3.1}, \ref{unthc} and \ref{cth} also consider quadratic BSDEs with concave generators. In addition, similar to Section 4 in \cite{Delbaen 2011}, we can also study the application of our Theorems \ref{unth3.1} and \ref{unthc} to quadratic PDEs, which is omitted in this paper due to the limitation of pages.

What we have seen from the above is the existence and uniqueness of unbounded solutions for quadratic BSDEs with generator $g$ being convex (concave) in $z$. The case of a non-convex generator was tackled in \cite{Fan-Hu-Tang 2020}, but more assumptions are required on the terminal value, that is, $|\xi|$ admits all exponential moments. When $|\xi|$ admits a certain exponential moment, the existence and uniqueness of unbounded solutions for quadratic BSDEs with non-convex generators are still unsolved, and this will be a subject of our future research.

%
%
%
%

\subsection*{Acknowledgment}
The authors would like to thank the anonymous referees for their careful reading
and some helpful suggestions.


\begin{thebibliography}{99}
		
\bibitem{Briand 2003}
P. Briand, B. Delyon, Y. Hu, E. Pardoux, L. Stoica, $L^p$ solutions of backward stochastic differential equations, {\it Stochastic Process. Appl.} {\bf 108} (2003) 109-129.

\bibitem{Briand-Hu 2006}
P. Briand, Y. Hu, BSDE with quadratic growth and unbounded terminal value, {\it Probab. Theory Relat. Fields} {\bf 136}  (2006) 604-618.

\bibitem{Briand-Hu 2008}
P. Briand, Y. Hu, Quadratic BSDEs with convex generators and unbounded terminal conditions, {\it Probab. Theory Relat. Fields} {\bf 141} (2008) 543-567.

\bibitem{Delbaen 2011}
F. Delbaen, Y. Hu, A. Richou, On the uniqueness of solutions to quadratic BSDEs with convex generators and unbounded terminal conditions, {\it Ann. Inst. Henri Poincar$\acute{e}$} {\bf 47} (2011) 559-574.

\bibitem{Delbaen 2015}
F. Delbaen, Y. Hu, A. Richou, On the uniqueness of solutions to quadratic BSDEs with convex generators and unbounded terminal conditions: the critical case, {\it Discrete Contin. Dyn. Syst.} {\bf 35} (2015) 5273-5283.

\bibitem{Fan 2016}
S. Fan, Bounded solutions, $L^p (p > 1)$ solutions and $L^1$ solutions for one-dimensional BSDEs under general assumptions, {\it Stochastic Process. Appl.} {\bf 126} (2016) 1511-1552.

\bibitem{Fan-Hu-Tang 2020}
S. Fan, Y. Hu, S. Tang, On the uniqueness of solutions to quadratic BSDEs with non-convex generators and unbounded terminal conditions, {\it C. R. Math. Acad. Sci. Paris} {\bf 358} (2020) 227-235.

\bibitem{Fan-Hu-Tang 2023}
S. Fan, Y. Hu, S. Tang, Multi-dimensional backward stochastic differential equations of diagonally quadratic generators: the general result, {\it J. Differ. Equations} {\bf 368} (2023) 105-140.

\bibitem{Fan-Jiang-Tian 2011}
S. Fan, L. Jiang, D. Tian, One-dimensional BSDEs with finite and infinite time horizons, {\it Stochastic Process. Appl.} {\bf 121} (2011) 427-440.


\bibitem{Hu 2005}
Y. Hu, P. Imkeller, M. Muller, Utility maximization in incomplete markets, {\it Ann. Appl. Probab.} {\bf 15} (2005) 1691-1712.

\bibitem{Hu-Tang-Wang 2022}
Y. Hu, S. Tang, F. Wang, Quadratic $G$-BSDEs with convex generators and unbounded terminal conditions, {\it Stochastic Process. Appl.} {\bf 153} (2022) 363-390.

\bibitem{Kobylanski 2000}
M. Kobylanski, Backward stochastic differential equations and partial differential equations with quadratic growth, {\it Ann. Probab.} {\bf 28} (2000) 558-602.

\bibitem{Luo-Fan 2018}
H. Luo, S. Fan, Bounded solutions for general time interval {BSDE}s
  with quadratic growth coefficients and stochastic conditions, Stoch. Dynam.
  {\bf 18} (2018) no. 1850034, 24pp.

\bibitem{Pardoux-Peng 1990}
E. Pardoux, S. Peng, Adapted solution of a backward stochastic differential equation, {\it Syst. Control Lett.} {\bf 14} (1990) 55-61.

\bibitem{Rouge 2000}
R. Rouge, N. El Karoui, Pricing via utility maximization and entropy, {\it Math. Finance} {\bf 10} (2000) 259-276.

\bibitem{Tian 2023}
D. Tian, Pricing Principle via Tsallis Relative Entropy in Incomplete Markets, {\it SIAM J. Financial Math.} {\bf 14} (2023) 250-278.

		
\end{thebibliography}
\end{document}